\documentclass[10pt,reqno]{amsart}
\usepackage{graphicx}
\usepackage{indentfirst,csquotes}

\usepackage{amssymb,amsthm,amsmath}
\usepackage{xcolor,paralist,hyperref,fancyhdr,etoolbox}

\theoremstyle{plain}

\newtheorem{theorem}{Theorem}[section]

\newtheorem{lemma}[theorem]{Lemma}
\newtheorem{proposition}[theorem]{Proposition}

\theoremstyle{remark}

\newtheorem{assumptions}[theorem]{Assumptions}

\newtheorem{definition}[theorem]{Definition}

\newtheorem{remark}[theorem]{Remark}

\hypersetup{ colorlinks=true, breaklinks=true, linkcolor=blue, citecolor = [RGB]{193, 37, 37} }

\newcounter{assumption-counter}%

\newcommand{\Prob}[2][]{\mathbb{P}{#1} \left[ #2 \right]}

\newcommand{\Esp}[2][]{\mathbb{E}{#1} \left[ #2 \right]}

\newcommand{\Tr}[1]{\operatorname{Tr} \left( #1 \right)}

\newcommand{\diag}[1]{\operatorname{diag} \left( #1 \right)}

\newcommand{\Id}[0]{\mathrm{Id}}

\newcommand{\Abs}[1]{\left\lvert #1 \right\rvert}

\newcommand{\Norm}[1]{\left\lVert #1 \right\rVert}

\newcommand{\Braces}[1]{\left\lbrace #1 \right\rbrace}

\newcommand{\Brackets}[1]{\left[ #1 \right]}

\newcommand{\Parentheses}[1]{\left( #1 \right)}

\newcommand{\Angles}[1]{\left\langle #1 \right\rangle}

\newcommand{\Ind}[1]{\mathbf{1}_{ #1 } }

\DeclareMathOperator{\diff}{d}

\DeclareMathOperator{\Diff}{D}

\DeclareMathOperator*{\argmin}{arg\,min}

\newcommand{\cadlag}{\emph{càdlàg}}

\newcommand{\EE}[0]{\mathbb{E}}
\newcommand{\FF}[0]{\mathbb{F}}

\newcommand{\HH}[0]{\mathbb{H}}

\newcommand{\NN}[0]{\mathbb{N}}
\newcommand{\PP}[0]{\mathbb{P}}

\newcommand{\RR}[0]{\mathbb{R}}

\newcommand{\Aa}[0]{\mathcal{A}}
\newcommand{\Bb}[0]{\mathcal{B}}

\newcommand{\Dd}[0]{\mathcal{D}}
\newcommand{\Ee}[0]{\mathcal{E}}
\newcommand{\Ff}[0]{\mathcal{F}}

\newcommand{\Hh}[0]{\mathcal{H}}

\newcommand{\Ll}[0]{\mathcal{L}}
\newcommand{\Mm}[0]{\mathcal{M}}

\newcommand{\Pp}[0]{\mathcal{P}}

\newcommand{\Rr}[0]{\mathcal{R}}
\newcommand{\Ss}[0]{\mathcal{S}}

\newcommand{\Ww}[0]{\mathcal{W}}
\newcommand{\Xx}[0]{\mathcal{X}}

\begin{document}


\title[Mean-Field Games with common Poissonian noise]{Mean-Field Games with common Poissonian noise: a maximum principle approach} 
\date{\today}
\thanks{This work was supported by the National Council of Science and Technology (CONACyT), scholarship number 863210, as part of the PhD Thesis of the second author.}


\author{Daniel~Hernández-Hernández}
\address{Centro de Investigación en Matemáticas A.C. Calle Jalisco s/n. 36240 Guanajuato, México}
\email{dher@cimat.mx}

\author{Joshué Helí~Ricalde-Guerrero}
\address{Centro de Investigación en Matemáticas A.C. Calle Jalisco s/n. 36240 Guanajuato, México}
\email{joshue.ricalde@cimat.mx}


\begin{abstract}
The theory of Mean-Field Games is interested in the behaviour of interacting particle systems in which the individual interaction between particles (players) decreases as the size of the population increases.  In recent years, it was introduced an interesting structure for this type of games,   assuming a correlated continuous source of randomness, which are called \textit{Mean-Field Games with Common Noise}. In this paper, we extend this concept and provide a precise definition of  a Mean-Field Game with Common Poissoninan Noise and its equilibrium. That is,  a common self-exciting Poissonian structure is considered within the dynamics of the population. Then, we address the problem of \textit{optimization for jump-diffusions with random environments} that lies within the definition of the MFG, and develop a stochastic version of the Pontryagin's Maximum Principle  to obtain a set of necessary conditions that an optimal control must satisfy.  Under additional convexity assumptions it is also shown that these conditions are also sufficient.


\end{abstract} 


\maketitle



\section{Introduction}
\label{Sec:Introdution}

Since their introduction at the beginning of the century in a series of pioneer works by Lasry \& Lions (\cite{lasry_jeux_2006}, \cite{lasry_jeux_2006-1}, \cite{lasry_mean_2007}) and, independently, by Caines, Huang \& Malhamé (\cite{minyi_huang_nash_2007}, \cite{huang_invariance_2007}, \cite{caines_large_2006}), \textit{Mean-Field Game} (MFG) theory has gained a great amount of appeal, becoming a relevant topic of interest that goes beyond only Game theory. This is due, in part,  to their flexibility in modelling interactions between (rational) agents, and the  advantages  they offer over  the so-called ``curse of dimensionality'' inherent to systems with a large number of participants. In addition to  their ability to capture both the \emph{macroscopic} (at a community-level) and \emph{microscopic} (at a representative individual-level) properties of a population.

MFGs serve as natural theoretical approximations for the study, control and optimization of very large dynamical systems. Heuristically,  a \emph{Mean-Field Equilibrium} (MFE) maintains the same notion of optimality as that of the original game. Formally, a MFE is the limit of the average state of a sequence of optimally controlled processes (in the sense of Nash equilibria), derived from symmetrical, large-population games \cite{cardaliaguet_master_2019}. However, the real utility and relevance of the theory comes from the fact that there is a similar result in the other direction: the policies obtained from a MFE form an $\epsilon$-Nash equilibrium when they are chosen by the players in the corresponding large-population game (\cite{minyi_huang_nash_2007}, \cite{carmona_probabilistic_2013}).  Clearly, this is an additional motivation to study MFGs derived from a large class of stochastic systems.

In order to further extend and broaden the scope of the usual MFG, Carmona, Delarue \& Lacker \cite{carmona_mean_2016} introduced an interesting structure in this type of games, assuming correlated randomness in a general environment; these are called {\it mean-field games with common noise}. Simply put, a MFG with Common Noise is one in which the stochastic dynamics of each particle is not only affected by an individual ({\it idiosyncratic}) continuous randomness source, but also by an external one that affects the whole system. As one can imagine, this model stands out within mean-field theory not only for its scope and practicality, but also for its complexity.
 
Historically, the study of MFGs with Common Noise has been centered around games in which the common noise process is  Gaussian, as it was  originally presented in \cite{carmona_mean_2016}. To the best of our knowledge, the most general result in this direction was obtained for a Gaussian white noise field $W^0 = \Parentheses{ W^0( \Lambda, B ), \Lambda \in \Bb(\Xi), B \in \Bb(\RR_+) }$; see the section of \emph{Notes and Complements} in \cite[Ch. 2]{carmona_probabilistic_2018-1}. However, while it is true that continuous (Gaussian) shock models provide a flexible set of tools for many situations, there are many others in which assuming a continuous interaction with the environment is unrealistic.
Take, for instance, the case of \emph{Systemic Risk} in economy and finance \cite{fouque_handbook_2013}. In the game theoretic formulation of the problem, one can interpret each particle as a \emph{bank} (\emph{agent}) within a \emph{banking (financial) market}. Then, the problem of finding the optimal investment strategy for an anonymous agent that is immersed in a tumultuous market can be translated into the MFG framework. This model has been previously  studied in the context of Mean-Field Games \cite{carmona_mean_2015}; however, the way the environment has been modeled (i.e. through continuous interactions) leaves out situations like ``black swans''; that is, uncommon or rare events which do not occur frequently, yet directly affect the entire market (such as catastrophes like the global pandemic of COVID-19), and are likely to have lasting effects thanks to the presence of a self-exciting component.

Jumps are by their very nature one of the best suited ways to model atypical and unexpected events. This fact   has not gone unnoticed, of course, as many recent studies on MFGs have indeed considered the inclusion of jumps into their models \cite{shen_maximum_2013}, \cite{hafayed_mean-field_2014}, \cite{bensoussan_mean-field-type_2020}, \cite{moon_linear-quadratic_2022}, \cite{xu__2021}, \cite{benazzoli__2019}, \cite{zhang_general_2018}, \cite{ledger_at_2021}.  Unfortunately, it is clear that this type of game has received less attention than its continuous counterpart.
Interestingly, an intrinsic characteristic   that these papers shared,   and which is not obvious to notice at first sight, is that none of them considers, implicitly or explicitly, the coupling of their corresponding dynamics through \emph{a ``big'' discontinuous process}. In other words, they assume that the propagation of chaos from their systems results in a \emph{McKean-Vlasov jump-diffusion} that is conditional to a Brownian motion filtration, which in turn can be traced  back to the  same models  discussed above. In fact, to the best of  our knowledge, the only work that consider a similar noise to the one we describe is \cite{barreiro-gomez_semiexplicit_2020}, where the common jump is given by an independent random measure independent from the rest of the elements in the system.

The election of jump processes depending on the environment is motivated by the need to obtain a sufficiently robust and  workable model that considers the effect of external inputs through periodic shocks. In many real-world applications, systems  are constantly affected by external factors, which do not necessarily evolve in continuous-time, and changes are observed after  certain positive random times, e.g. new arrivals to  queue systems or claims to assurance companies. In these cases, it is reasonable to consider  that periodic shocks are directly affected by the existing environment (arrivals diminish if there is a congestion before the entrance to the queue, and claims increase if the locality is exposed to natural disasters). By allowing the structure of the jumps process to have a relationship with the environment it is possible to incorporate exogenous information. This relationship is established in terms of the intensity of the jumps process, covering a large number of problems, while maintaining  meaningful results.


With this in mind, the aim of this paper is to answer the following question: \textit{Assuming that a Mean-Field Game can be derived from the limit of a large symmetric stochastic differential game such that its players are all affected by a common discontinuous component with Poissonian structure (self-exciting or marked), when and how are we able to solve it?}
In more technical terms, we start from our previous results on the existence of a solution to Forward-Backward SDEs immersed in randomized environments and their propagation of chaos (\cite{hernandez-hernandez_coupled_2023} and \cite{hernandez-hernandez_conditional_2023}, respectively) and introduce Mean-Field Games in which the coupling between players is reinforced by a common discontinuous process. Then,  variational methods are used to deduce \textit{necessary conditions} that an equilibrium must satisfy (namely, an appropriate version of \textit{Pontryagin's Maximum Principle}) such that, when  additional convexity and compactness conditions are included, these also become  \textit{sufficient}.

The paper is structured as follows. The notation and  main elements of the model are introduced in Section \ref{Sec:Basic-Notation}. Section \ref{Sec:Heur} is dedicated to motivate the main problem, starting with the description of a game with a finite number of players, and analyzing the limit when the population size grows to infinity. With this motivation at hand,   we proceed to describe precisely the Mean-Field Game we are interested in,  in particular its equilibrium, with special emphasis on the dynamics describing its evolution. The maximum principle is stated  in Section \ref{Sec:PMP}, providing necessary and sufficient conditions for the equilibrium strategy. The solution to the Mean-Field Game is presented in Section \ref{Sec:Sol MFG} in terms of the solution of the Conditional McKean-Vlasov FBSDE. Finally, two important examples are discussed in Section \ref{Sec:EyA}, including the LQ model, and the mean-variance portfolio selection within regime switching.

\section{Basic notation and preliminaries}
\label{Sec:Basic-Notation}

Given $x \in \RR^d$, we write $x^{-i} = \Parentheses{ x^1, \ldots, x^{i-1}, x^{i+1}, \ldots, x^d }$, $i=1, \ldots, d$, to represent the vector $x$ without the $i$-th coordinate. The usual inner product in $\RR^d$ is denoted by $(x\cdot y) = x^\top y$.  If $G$ is a given full rank matrix on $\RR^{n \times d}$, we write $(x\cdot y)_G := (Gx)^\top y$ for  $x \in \RR^d$ and  $y \in \RR^n$. We also use $c_G \in \RR$ to denote the positive proportionality constant $( x \cdot Gx)_G = c_G (x\cdot x)$. For matrices, we use the Frobenius inner product, i.e. $(A \cdot B)_\mathrm{Fr} = \Tr{A B^\top}$ for $A,B$ real-valued matrices such that the product is well-defined, with the corresponding norm given by $\Norm{A}_\mathrm{Fr}^2 = \Tr{ A A^\top }$.

Given the time interval $[0,T]$, with $T>0$ fixed, and a Polish space $(\Xx,d_\Xx)$, we denote by 
\begin{align*}
    \mathrm{C}([0,T];\Xx)
    :=&
    \Braces{X:[0,T] \to \Xx \ \middle\vert \ X\mbox{ has continuous trajectories }},
    \\
    \mathrm{D}([0,T];\Xx)
    :=&
    \Braces{X:[0,T] \to \Xx \ \middle\vert \ X\mbox{ has \cadlag\ trajectories }}.
\end{align*}
We also denote the space of measures $\Mm$ and $\Mm_c^*$ as
\begin{align*}
    \Mm(\Xx) :=& \Braces{ \nu \middle\vert \nu\mbox{ is a  finite measure over }(\Xx,\Bb(\Xx)) },
    \\
    \Mm_c^*(\Xx) :=& 
        \Braces{ \nu \in \Mm(\Xx) \middle\vert \nu\mbox{ is a simple, counting measure over }(\Xx,\Bb(\Xx)) }.
\end{align*}
Product spaces will always be endowed with the product $\sigma-$field. Lastly, we use $\Dd([0,T];\Xx) = \Bb\Parentheses{ \mathrm{D}([0,T];\Xx) }$ for the Borel $\sigma$-algebra with respect to the Skorohod topology \cite{billingsley_convergence_2013}.

\subsubsection*{Wasserstein spaces.} Denote by $\Pp(\Xx)$ the set of Borel probability measures on $\Xx$ and, for $p \in [1,\infty)$, define the set of probability measures with finite moments of order $p$ as
\begin{align*}
    \Pp_p(\Xx)
    := 
    \Bigg\lbrace 
        \nu \in \Pp(\Xx) 
        \ \Bigg\vert\ 
        \int_\Xx d_\Xx(x_0,x)^p \nu( \diff x ) < +\infty \  \mbox{for  any } x_0 \in \Xx
    \Bigg\rbrace.
\end{align*}
Given $\nu^1,\; \nu^2\in \Pp(\Xx) $, define $\Pi(\nu^1,\nu^2) \subset \Pp(\Xx \times \Xx)$ as the set of probability measures over the product space $\Xx \times \Xx$ with marginals $\nu^1$ and $\nu^2$, respectively. The $p$-th Wasserstein distance between two probability measures $\nu^1, \nu^2 \in \Pp_p(\Xx)$ is then defined as
\begin{align}
    \label{Eq:W_p^p}
    W_p (\nu^1, \nu^2)^p
    :=
    \inf \Braces{ 
        \int_{ \Xx \times \Xx } d_\Xx( x_1, x_2 )^p \pi( \diff x_1, \diff x_2 )
        \ \middle\vert\ 
        \pi \in \Pi( \nu^1, \nu^2 )
    }.
\end{align}
When $p=2$, we will be using repeatedly the following inequality:
\begin{align}
    \label{Eq:W_2-Inequality}
    W_2 (\nu^1, \nu^2)^2 \leq \Esp{ d_\Xx(Y^1,Y^2)^2 },
\end{align}
where $Y^1,Y^2$ are independent $\Xx$-valued random variables defined on the same probability space $(\Omega,\Ff,\PP)$ and $\nu^i = \PP \circ (Y^{i})^{-1}$ for $i=1,2$.

Any mapping $\Psi:\Omega \times [0,T] \times \Pp_2(\RR^d) \times \RR^d \to \RR^n$ is denoted by $\Psi_t(\nu,x)(\omega)$ for each  $(\nu,x) \in \Pp_2(\RR^d) \times \RR^d$ and $(\omega, t) \in \Omega \times [0,T]$, and by $\Psi_t^\nu(x)(\omega)$ when there is no ambiguity regarding the measure term. Moreover, for  each $(\nu,x),(\nu',x')\in \Pp_2(\RR^d) \times \RR^d$, their difference is denoted by $ \delta \Psi_t(\nu,x)(\omega)= \Psi\Parentheses{ \omega, t, \nu, x } - \Psi\Parentheses{ \omega, t, \nu', x' }$. When the difference of the mapping $\Psi$ only occurs on the $x$ term (respectively, only on the $i$-th component of $x$), we write 
\begin{align*}
    \delta \Psi_t^{\nu}\Parentheses{ x }(\omega) 
    :&= 
    \Psi\Parentheses{ \omega, t, \nu, x } - \Psi\Parentheses{ \omega, t, \nu, x' }
    \\
    (\mbox{resp. }\delta \Psi_t^{\nu}\Parentheses{ x^i, (x^i)' \middle| x^{-i} }(\omega)
    :&=
    \Psi\Parentheses{ \omega, t, \nu, x^i, x^{-i} } - \Psi\Parentheses{ \omega, t, \nu, (x^i)', x^{-i} }).
\end{align*}
This notation is also used to represent the difference of any two elements of an Euclidean space $x,x'$  as $  \delta x = x - x'$. If in addition $\Psi$ is twice continuously differentiable on the $x$ component, we denote by $\Diff\Psi$ (resp. $\Diff^2\Psi$) the Jacobian (resp. Hessian) of $\Psi$ with respect to $x$. Moreover, in general, we define the sets
\begin{align*}
    \mathrm{C}^p(\RR^d;\RR^n)
    :=&
    \Braces{ 
        \Psi \in \mathrm{C}(\RR^d;\RR^n)
        \ \vert\ 
        \Psi\mbox{ is }p\mbox{-times continuously  differentiable}},
    \\
    \mathrm{C}^p_b(\RR^d;\RR^n)
    :=&
    \Braces{ 
        \Psi \in \mathrm{C}^p(\RR^d;\RR^n)
        \ \middle\vert\ 
        \Diff^{k} \Psi \mbox{ is bounded for }k=0,1,\ldots,p
    },
\end{align*}
for $p \in \NN$, where $\Diff^{k}\Psi$ denotes the $k$-th differential of $\Psi$ with respect to $x$. In each  case, we omit the last term in the parentheses when there is no ambiguity regarding the image space.

\subsubsection*{Spaces for stochastic processes.} Given   a stochastic process $\Theta$, we will some times need to introduce other filtrations associated with $\Theta$. The natural filtration is denoted by $\FF^{ \mathrm{nat},\Theta }$, while the $\PP$-null complete right-continuous augmentation is denoted by $\FF^{ \Theta }$. In general, we say that a filtration $\FF$ satisfies the usual conditions if it is right-continuous and complete.

We write $\Hh^2(\RR^d)$ (or simply $\Hh^2$ when there is no ambiguity) for the Hilbert space of $\RR^d$-valued, square-integrable martingales (with respect to $\FF$), i.e.
\begin{align*}
    \Hh^2(\RR^d)
    :=
    \Bigg\lbrace 
		M:\Omega \times &[0,T] \to \RR^d
		\ \Bigg\vert\ 
        M\mbox{ is an }\FF\mbox{-martingale with } \sup_{ t \in [0,T] }\Esp{ \Abs{ X_t }^2 }<\infty
    \Bigg\rbrace.
\end{align*}
This space is endowed with  the norm $\Norm{\ \cdot\ }_{\Hh^2(\RR^d)}$, induced by  the inner product $\Parentheses{ M \cdot M' }_{\Hh^2}=\Esp{ M_T^\top M'_T }$, defined  by $\Norm{ M }_{ \Hh^2 }:=\Esp{ \Abs{ M_T }^2 }^{ \frac{1}{2} }$.  For any subset $B \subset \Hh^2$, $B^\perp \subset \Hh^2$ denotes the orthogonal set of $B$, i.e. the set of elements $M\in \Hh^2$ being  orthogonal to some martingale in $B$. The space of well-posed integrands with respect to $M$  is denoted by $\HH^2(M)$, and it is endowed with the norm of predictable quadratic variation $\Norm{\ \cdot\ }_{\HH^2(M)}$. If $M'$ is another square integrable martingale, we write $\HH^2(M,M') = \HH^2(M) \cap \HH^2(M')$. We denote by $\Ll^2(M)$ the space of martingale-driven, Itô stochastic integrals of processes in $\HH^2(M)$ with respect to $M$, i.e.  $\Ll^2(M):=\Braces{ \int H \diff M\ :\ H \in \HH^2(M) },$ equipped with the martingale norm \cite[Ch. IV.2]{protter_stochastic_2005}.  We use the notation $[H,H']$ to denote the \textit{quadratic covariation} between semimartingales $H$ and $H'$, and $\Angles{H,H'}$ to denote its compensator when it exists; this process is also known as \textit{conditional} or \textit{predictable quadratic covariation} \cite[Ch. III.5]{protter_stochastic_2005}. In particular,  we write $[H]=[H,H]$ and $\Angles{H} = \Angles{H,H}$. The space of adapted, $\Xx$-valued \cadlag\ processes with finite supremum norm is denoted by $(\Ss^2(\Xx), \Norm{\ \cdot\ }_\ast$). In the $d$-dimensional case,
\begin{align*}
    &\Norm{ M }_*^2 = \max_{1 \leq i \leq d} \Norm{ M^{(i)} }_*^2 
        = \max_{1 \leq i \leq n} \Esp{ \sup_{[0,T]} \Abs{ M^{(i)}_t }^2 },
\end{align*}
for all $M = (M^{(1)},\ldots,M^{(d)}) \in \Ss^2(\Xx)$. Moreover, $\Parentheses{ \Ss^2(\Xx), \Norm{\ \cdot\ }_* }$ is a Banach space \cite[Ch. V.2]{protter_stochastic_2005}. We use the symbol $\FF' \hookrightarrow \FF$ whenever $\FF'$ and $\FF$ are filtrations on the $\sigma$-algebra $\Ff$ such that $\FF' \subset \FF$ and  $\FF'$-martingales are also $\FF$-martingales. In particular, an $\FF$-adapted \cadlag\ process $\Theta = \Parentheses{ \Theta_t, 0 \leq t \leq T }$ with values in a Polish space is said to be \emph{compatible with $\FF$ (under $\PP$)} if $\FF^\Theta\hookrightarrow\FF$.

\subsubsection*{Admissibility}

By their very nature, the interactions described by the mean-field model produce  measured-valued processes that act as random environments for the particles themselves. Addressing this problem requires answering the following questions: How to deal at the same time with a discontinuous process and a random (possibly exogenous) environment, both defined in a general probability space, and which feedback on each other at all times? This problem is precisely what gives rise to the notion of \textit{admissibility}.

The idea of an admissible set-up is not new to MFGs (see \cite[Ch. 1]{carmona_probabilistic_2018-1}). However, unlike the case of a common Gaussian noise, in this paper   the inclusion of the predictable quadratic variation of the common noise   is included \cite{hernandez-hernandez_coupled_2023}. Heuristically, \textit{admissibility is a way to moderate the effects of the environment through the compensator} and, as a consequence, it determines how the underlying probability space is defined, how the interactions and the convergence are obtained, and also how the optimal control is chosen. 


Formally, a random environment $\mu$ refers to a measurable map defined on some probability space $(\Omega, \Ff, \PP)$ into 
\begin{align*}
    \Parentheses{ 
		\mathrm{D} \Parentheses{ [0,T]; \Pp_2(\RR^d) }, 
		\Dd \Parentheses{ [0,T]; \Pp_2(\RR^d) } 
    },
\end{align*}
with $\FF^\mu \subset \FF$, such that there exists an $\Ff_0^{\mu}$-measurable variable $X_0$ and a probability measure $\nu_0 \in \Pp_2(\RR^d)$, representing the \emph{initial distribution of $\mu$}, that satisifies 
\begin{align}
	\label{Eq:Conditions-on-Random-Environment-I}
	\Ll(X_0) &= \nu_0,
	&
	\Prob{  \mu_0 = \nu_0} &= 1,
\end{align}
and, for each $\nu\in\Pp_2(\RR^d)$,
\begin{align}
	\label{Eq:Conditions-on-Random-Environment-II}
	\Esp{ \sup_{t\in [0,T]} W_2(\mu_t,\nu)^2 } < \infty.
\end{align}

\begin{remark}
In the context of MFGs, $\mu$ concentrates along the time the states of the entire set of interacting particles into a single population statistic.
\end{remark}

Since we are interested in \emph{marked Poisson processes} with random intensity and marks on a measurable space $(\mathbf{R},\Rr)$ such that the information from the environment influences their behavior, we will be assuming first that $\FF^{\mu} \subset \FF$ for some general filtration $\FF$, and that $\lambda : \Omega \times [0,T] \times \mathbf{R} \times \Pp_2(\RR^d) \to \RR^l_+$ is an \textit{$\FF$-admissible intensity candidate}, which means that $\lambda$ satisfies the following conditions.
\begin{enumerate}
    
    \item
        $\lambda_t^{(i)}(\cdot,r,\nu) \perp \lambda_t^{(j)}(\cdot,r,\nu)$ for $i \neq j$, and each $(t,r,\nu) \in [0,T] \times \mathbf{R} \times \Pp_2(\RR^d)$;
    
    \item 
        $\lambda$ is $\FF$-predictable;
    
    \item
        for all $0 \leq a < b \leq T$,
        \begin{align*}
            &\int_a^b \int_{\mathbf{R}} \lambda_t^{(i)}(r,\nu) Q^{(i)}(\diff r) \diff t < \infty
            &
            &\PP-\mbox{a.s.},
    \end{align*}
    where $Q^{(1)},\ldots,Q^{(l)}$ are  finite measures on $(\mathbf{R},\Rr)$.
    
\end{enumerate}
Then, we say that the random kernel $K^\lambda : \Omega \times [0,T] \times \Rr \to \RR_+^l$, defined entry-wise as
\begin{align}
    \label{KerK}
    &K^{\lambda(i)}_t(\omega, C) := \int_{C} \lambda^{(i)}_t(\omega,r,\mu_{t-}(\omega)) Q^{(i)}(\diff r),
    &
    &i = 1, \ldots, l,
\end{align}
for $(\omega,t) \in \Omega \times [0,T]$ and  $C \in \Rr$, is an \emph{admissible $\FF$-intensity kernel}.

\begin{remark}
\label{Remark:Form-of-lambda}
Using  a slight abuse of notation, observe that a particular case of admissible candidates correspond to processes of the form 
\begin{align*}
    &\lambda_t(\omega,r) = \lambda(r,\mu_{t-}(\omega)) 
    &
    &\forall (\omega,t,r) \in  \Omega \times [0,T] \times \mathbf{R},
\end{align*}
for some deterministic function $\lambda$; i.e. candidates where the randomness only occurs through the environment at the previous instant. In this manuscript we focus  on this kind of candidates, since they fit well with the random jump model we are interested in.
\end{remark}

Consequently, an \textit{admissible noise process} is defined as a marked Poisson process $(N^\lambda,\xi)$ in the time interval $[0,T]$ with marks on $\mathbf{R}$, and $\FF$-intensity kernel $K^\lambda$. This means that $N^\lambda$ is a Poisson process with stochastic intensity $K^\lambda_t(\omega,\mathbf{R}) = \int_\mathbf{R} \lambda_t(\omega,r)Q(\diff r)$ and jumping times $\tau = \Braces{\tau_n, n \geq 1 }$:
\begin{align*}
    \tau_{n + 1}(\omega)
    = 
    \inf\Braces{ t \geq \tau_{n}(\omega) ; \int_0^{t} K^\lambda_s(\omega,\mathbf{R}) \diff s \geq {n+1}}, \;\;\PP-\mbox{a.s.}.
\end{align*}
The second component, $\xi = \Braces{ \xi_i, i \geq 1 }$, is a sequence of $\mathbf{R}$-valued random variables with transition probabilities
\begin{align*}
    \Prob{ \xi_n \in C \middle\vert \sigma\Braces{ (\tau_i, \xi_i), i = 1, \ldots, n-1 } \vee \sigma\Braces{ \tau_n } }(\omega)
    =
    \frac{ K^\lambda_{\tau_n}(\omega, C ) }{ K^\lambda_{\tau_n}(\omega, \mathbf{R} ) }.
\end{align*}
Note that $(N^\lambda,\xi)$ can also be characterised by its so-called \textit{lifting} \cite[Ch. 5]{bremaud_point_2020}; that is, the random measure $\eta^\lambda$ defined as
\begin{align}
    \label{Eq:Lifting}
    &\eta^\lambda(B) = \sum_{i \geq 1} \delta_{(\tau_i,\xi_i)}(B),
    &
    & B \in \Bb([0,T]) \otimes \Rr,
\end{align}
with predictable compensator $K^\lambda_t(\omega, \diff r) \diff t$.


\begin{remark}
In general, when we refer to $(N^\lambda,\xi)$ as an $l$-dimensional marked Poisson process, we refer to the family
\begin{align*}
    (N^\lambda,\xi) = \Braces{ (N^{\lambda(j)},\xi^{(j)}), j = 1, \ldots, l },
\end{align*}
where each $(N^{\lambda(j)},\xi^{(j)})$ is an independent Poisson process with an $\FF$-intensity kernel $K^{\lambda(j)}$ from $(\Omega \times [0,T], \Ff \otimes \Bb([0,T]))$ into $(\mathbf{R},\Rr)$. Similarly, we write $K^{\lambda}$ (resp. $\eta^\lambda$) for the vector with $j$-th entry $K^{\lambda(j)}$ (resp. $\eta^{\lambda(j)}$); see \eqref{KerK}-\eqref{Eq:Lifting}. 
\end{remark}

We can now describe an \emph{admissible Set-up} as an space endowed with both an \emph{idiosyncratic} and  \emph{admissible noises}.

\begin{definition}[\cite{hernandez-hernandez_coupled_2023}]
\label{Def:Admissible-Set-Up}
A tuple $(\Omega,\Ff,\FF,\PP,X_0,\mu,K^\lambda,(N^\lambda,\xi),W)$, or simply $( X_0, \mu, K^\lambda, ( N^\lambda, \xi )$, $W )$, is said to be an \emph{admissible set-up}  if
\begin{enumerate}
    \item
	$\mu$ is a random environment on $(\Omega, \Ff, \PP)$;  
    \item
	$X_0 \in L^2( \Omega, \Ff_0, \PP; \RR^d )$;
		
    \item
        $K^\lambda$ is a non-negative, $\FF$-predictable, locally integrable kernel from $\big( \Omega \times [0,T],\Ff \otimes \Bb([0,T]) \big)$ into the measurable space $(\mathbf{R},\Rr)$;
		
    \item
        $(N^\lambda,\xi)$ is an $l$-dimensional marked Poisson process, adapted to $\FF$, with $\FF$-stochastic intensity kernel $K^\lambda$ and lifting $\eta^\lambda$;

    \item
        A sequence  $W=\{W^i\}$ of independent $k$-dimensional $\FF$-Wiener processes, referred as an \emph{idiosyncratic noise};
		
    \item
        $(X_0, \mu, \eta^\lambda)$ is independent of $W$ (under $\PP$);
		
    \item
	Any $\FF^{(X_0, \mu, \eta^\lambda, W)}$-martingale is also an $\FF$-martingale (under $\PP$).
\end{enumerate}
\end{definition}

\begin{remark} 
\label{Remark:Intensity-Candidate}
\noindent
\begin{enumerate}
    
    \item
        Condition 7 above is commonly known in the literature as \emph{$\Hh$-hypothesis}, meaning that $\FF^{(X_0, \mu, \eta^\lambda, W)}\hookrightarrow\FF$, or simply that ${(X_0, \mu, \eta^\lambda, W)}$ is \emph{compatible} with $\FF$; see \cite{aksamit_enlargement_2017}.
    
    \item 
        Since we are only considering kernels $K^\lambda$ that are obtained through an intensity candidate $\lambda$, we write $(X_0,\mu,\lambda,(N^\lambda,\xi),W)$ when referring to the admissible set-up $( X_0, \mu$, $K^\lambda, ( N^\lambda, \xi ), W)$.
    
\end{enumerate}
\end{remark}

In order to describe properly solutions to the Hamiltonian system of FBSDEs (appearing below) in a random environment, we need to define the following space:
\begin{align}
    \label{Eq:Frac-H-lambda}
    \mathfrak{H}^\lambda
    :=
    \Big\lbrace
        U : \Omega \times [0,T] \times \mathbf{R} \to \RR^{n \times l}
		\Big\vert\ &
		U \text{ is } \FF \otimes \Rr-\text{predictable, with } 
        \\ \nonumber
		&\Norm{ U_t(\omega,\cdot) }_{\lambda_t ( \omega )} < \infty,
        \\ \nonumber
		&\diff \PP \otimes \diff t\mathrm{-a.e. }(\omega, t) \in \Omega \times [0,T]
    \Big\rbrace,
\end{align}
where $\Norm{\ \cdot\ }_{\lambda_t (\omega)}^2 = \Parentheses{ \cdot, \cdot }_{\lambda_t (\omega)}$ is a random norm defined for matrix-valued functions, induced by the inner product
\begin{align}
    \label{Eq:lambda-inner-product}
    &\Parentheses{ u \cdot v }_{\lambda_t (\omega)} 
    := 
    \int_\mathbf{R} \Tr{ u(r) \diag{K^\lambda_t(\omega,\diff r)} v(r)^\top },
\end{align}
for $\diff \PP \otimes \diff t$-a.a. $(\omega, t) \in \Omega \times [0,T]$. This is due to the fact that, for each $U \in \mathfrak{H}^\lambda$, if
\begin{align}
    \label{Eq:Finite-Predictable-Quadratic-variation-N^lambda}
    &\int_0^T \Norm{U_t(\cdot)}_{\lambda_t}^2 \diff t < \infty
    &
    &\PP-\mbox{a.s.},
\end{align}
then the integral process defined below is well-defined for $\diff \PP \otimes \diff t$-almost all $(\omega,t) \in \Omega \times [0,T]$:
\begin{align}
    \label{Eq:Integral-wrt-eta^lambda}
    \int_{[0,t] \times \mathbf{R} } U_s(\omega,r) {\eta}^\lambda &(\omega, \diff s, \diff r)
    :=
    \sum_{j=1}^l \sum_{i \geq 1} U^{(\cdot,j)}_{ \tau_i^{(j)} } (\omega,\xi_i^{(j)}(\omega)) \Ind{ \{\tau_i^{(j)} \leq t\} },
\end{align}
where $\tau^{(j)}$ (resp. $\xi^{(j)}$) is the sequence of jumping times (marks) of the $j$-th member of $(N^\lambda,\xi)$, and $U^{(\cdot,j)}$ denotes the $j$-th column of $U$. Moreover, if the random variable in \eqref{Eq:Finite-Predictable-Quadratic-variation-N^lambda} is also integrable with respect to $\PP$, then the integrals with respect to the compensated martingale measure
\begin{align}
    \label{Eq:tilde-eta^lambda}
    \tilde{\eta}^\lambda(\omega, \diff t, \diff r)
    :=
    {\eta}^\lambda(\omega, \diff t, \diff r)
    -
    K_t^{\lambda}(\omega, \diff r) \diff t
\end{align}
are also  square-integrable martingales, see \cite[Thm. 5.1.33-34]{bremaud_point_2020}. 

We refer to \eqref{Eq:Integral-wrt-eta^lambda} as the \emph{integral of $U$ with respect to $(N^\lambda,\xi)$} --equivalently, the \emph{integral of $U$ with respect to $\eta^\lambda$}--. Following the notation used for well-posed integrals with respect to square-integrable martingales, we write
\begin{align*}
    \HH^2(\tilde{\eta}^\lambda)
    :=
    \Braces{
        U \in \mathfrak{H}^\lambda
        \  \middle\vert\ 
        \Esp{ \int_0^T \Norm{U_t(\cdot)}_{\lambda_t}^2 } < \infty
    }
\end{align*}
for the space of well-posed integrands with respect to $(N^\lambda,\xi)$, and 
\begin{align*}
    \Ll^2(\tilde{\eta}^\lambda)
    :=
    \Braces{
        \Parentheses{
            \int_{[0,t] \times \mathbf{R} } U(s,r) \tilde{\eta}^\lambda (\diff s, \diff r)
        }_{t \in [0,T]}
        \ \middle\vert\ 
        U \in \HH^2(\tilde{\eta}^\lambda)
    }
\end{align*}
for the space of martingales that can be expressed as compensated integrals with respect to $(N^\lambda,\xi)$.


\section{Heuristic derivation of the Mean-Field Game}
\label{Sec:Heur}

Let $A$ be a fixed non-empty convex and compact subset of a finite-dimensional Euclidean space, which is referred hereafter as the \emph{set of admissible actions}. Given $n \in \NN$ (fixed) representing the number of players in the game, an element $a^{n} \in \mathrm{D}([0,T]; A^n )$ is called a \emph{strategy profile}.  Throughout we focus on \emph{autonomous} or \emph{time-invariant} systems where the time dependency occurs through the state of the system and the environment. Under the prevailing assumptions established in \cite{hernandez-hernandez_conditional_2023} and adapted to the controlled model, let $X^{a^n} \equiv X^{n} = \Braces{ X^{i,n}, 1 \leq i \leq n }$ be the unique solution process to
\begin{align}
    \label{Eq:n-player-dynamics}
    X^{i,n}_t
    =&
    X^{i,n}_0
    +
    \int_0^t b^{i,n} \Parentheses{ X^{n}_{s-}, a^{i,n}_{s-} } \diff s
    +
    \int_0^t \sigma^{i,n} \Parentheses{ X^{n}_{s-}, a^{i,n}_{s-} } \diff W_s^i 
    \\ \nonumber
    &+
    \int_{[0,t] \times \mathbf{R} } \gamma^{i,n} \Parentheses{ r, X^{n}_{s-}, a^{i,n}_{s-} } 
        \tilde{\eta}^{n} (\diff s, \diff r),
\end{align}
where $W^1,\ldots,W^n$ are independent $k$-dimensional Brownian motions, $\eta^{n}$ is the lifting of a marked Poisson process $(N^{n},\xi^{n})$ with transition kernel $\lambda^n( r, X^n_{t-} ) Q(\diff r)$, and $\tilde{\eta}^{n}$ is the corresponding compensated measure. We refer to $X^{a^n}$ as the stochastic differential system controlled by $a^n$.

In the context of differential games, each $X^{i,n}$ represents the \emph{state of the system controlled by the $i$-th player}, with the flow $a^{i,n}$ being their corresponding \emph{strategy}. The \emph{set of admissible strategies} for each player is therefore defined as
\begin{align*}
    \mathcal{A}
    :=
    \Braces{ 
        \alpha : \Omega \times [0,T] \to A
        \ \middle\vert \ 
        \alpha\mbox{ is \cadlag\ and }\FF-\mbox{adapted}
    },
\end{align*}
with each pair $(X^\alpha,\alpha)$ being referred to as an $\alpha$-controlled pair.  The \emph{running cost} and \emph{terminal cost} for the $i$-th player, $i = 1, \ldots, n$ are defined as a given pair of measurable functions
\begin{align*}
    &f^{i,n}: (\RR^d)^n \times A \to \RR,
    &
    &g^{i,n}: (\RR^d)^n \to \RR,
\end{align*}
respectively; consequently, for a given vector of strategies $\alpha^n$ the \emph{pay-off for the $i$-th player} is defined as the mapping
\begin{align}
    \label{Eq:n-player-payoff}
    \Ss^2(\RR^d)^n \times \Aa \ni (X^{n},\alpha^{i,n})
    \longmapsto &
    \Esp{ \int_0^T f^{i,n}( X_{t-}^n, \alpha^{i,n}_{t-} ) \diff t + g^{i,n}( X^{n}_{T-} ) }
    \\ \nonumber
    &=:
    J^{i,n}( X^{n}, \alpha^{i,n} ).
\end{align}

A \textit{Nash equlibrium for the $n-$player game} is defined as the collection of controlled pairs $(\widehat{X}^n,\widehat{\alpha}^n) = \{ (\widehat{X}^{i,n},\widehat{\alpha}^{i,n}), 1 \leq i \leq n \}$ such that the following condition is fulfilled for each player ($i=1,\ldots, n$):
\begin{align*}
    &J^{i,n}(\widehat{X}^{n},\widehat{\alpha}^{i,n}) 
        \leq J^{i,n}(
            \widehat{X}^{1,n},\ldots,\widehat{X}^{i-1,n},
            X,
            \widehat{X}^{i+1,n},\ldots,\widehat{X}^{n,n},
        \alpha),
\end{align*}
for all $(X,\alpha) \in \Ss^2(\RR^d) \times \Aa$, where $X$ is the  solution of the system (\ref{Eq:n-player-dynamics}) when an arbitrary strategy $\alpha$ is applied for player $i$. In other words: \textit{when the game is at an equilibrium, there is no reason for any individual player $i$ to unilaterally change their strategy from $\widehat{\alpha}^{i,n}$ to any other strategy $\alpha \in \Aa$} \cite{pontryagin_mathematical_1986}, \cite{yong_stochastic_1999}, \cite{boel_optimal_1977}, \cite{oksendal_applied_2019}. Furthermore, it is a known fact that under suitable Lipschitz and convexity assumptions on the coefficients, there exists a Nash equilibrium for the $n$ player game \eqref{Eq:n-player-dynamics}-\eqref{Eq:n-player-payoff}, and that it is possible to find a collection of deterministic measurable functions 
\begin{align}
    \label{Eq:Optimal-policy-I}
    &\widehat{\mathbf{a}}^{i,n} : [0,T] \times (\RR^{d})^n \to A,
    &
    &i=1,\ldots,n,
\end{align}
each one called an \emph{optimal policy}, such that the strategy profile given by
\begin{align}
    \label{Eq:Optimal-policy-II}
    &\widehat{\alpha}^{i,n}_t \equiv \widehat{\mathbf{a}}^{i,n}( t, \widehat{X}^{n}_t ),
    &
    &\forall t \in [0,T]
\end{align}
for each $i=1,\ldots,n$, is a Nash equilibrium for the game \cite{bensoussan_stochastic_1983}, \cite{peng_general_1990}, \cite{kabanov_pontryagin_1997}.


\begin{remark}
The above assumptions about the existence of a Nash equilibrium and its optimal policy  are quite strong and deserve their own study (see \textit{Markovian control policies}  in \cite{fleming_controlled_2006} and the references there-in).  However, it is not our intention to analyze this problem for the moment, and we only describe it as a means to motivate the idea of Mean-Field equilibria as a thermodynamic limit.
\end{remark}

Previous definitions are standard from game theory; nevertheless, the key aspect that differenciates a Mean-Field Game from the above is the idea that each player is \textit{virtually playing the same game against the rest}, and that \textit{none of them individually can directly affect the others except from their contribution to the average of the system}. Furthermore, the individual contribution from an individual player becomes negligible (and thus, the player becomes indistinguishable from the rest) as the size of the population grows. Mathematically, this means that the coefficients of the game \textit{are the same for every player} ($b^{i,n}=b^n$ for all $i$, and similarly for $\sigma^{i,n},\gamma^{i,n}, \ldots$), and that they satisfy the following notion of \textit{symmetry}: Let $x^n = ( x^{1,n}, \ldots, x^{n,n} ) \in (\RR^d)^n$ and $a \in A$; then, $b^n$ is permutation-invariant in $x^n$:
\begin{align}
    \label{Eq:Symmetric-coefficient}
    b^n( x^n, a ) = b^n( x^{ \pi(1),n }, \ldots, x^{ \pi(n),n }, a )
\end{align}
where $\pi$ is any arbitrary permutation on the indices $\Braces{ 1, \ldots, n }$. This symmetry is precisely what allows us to move from a \textit{finite population} regime to the \textit{empirical} one presented in \cite{hernandez-hernandez_conditional_2023}.

First, observe that under these symmetry assumptions, if $\widehat{\alpha}^n$ is a Nash equilibrium for the game \eqref{Eq:n-player-dynamics}-\eqref{Eq:n-player-payoff}, then the processes $\{ \widehat{\alpha}^{1,n},\ldots,\widehat{\alpha}^{n,n} \}$ are exchangeable. Additionally, if there exist optimal policies $\widehat{\mathbf{a}}^{i,n}$ as in \eqref{Eq:Optimal-policy-I}-\eqref{Eq:Optimal-policy-II}, then $\widehat{\mathbf{a}}^{i,n}$ is also symmetric in the sense of \eqref{Eq:Symmetric-coefficient}, with $\widehat{\mathbf{a}}^{n} = \widehat{\mathbf{a}}^{i,n}$  for all $i = 1, \ldots, n$.

Second, given    an arbitrary compact  subset $B \subset \RR^d$,  let $\{ \Psi^n:(\RR^d)^n \to \RR, n \geq 1 \}$ be a sequence of smooth symmetrical functions, in the sense of \eqref{Eq:Symmetric-coefficient}, such that
\begin{align*}
    &\Norm{ \Psi^n }_{L^\infty(B)} \leq C
    &
    &\mbox{and}
    &
    &\Norm{ \Diff \Psi^n }_{L^\infty(B)} \leq \frac{C}{n},\ \forall n,
\end{align*}
for some constant $C > 0$. Then, there exists a subsequence $\Braces{ \Psi^{ n_k }, k \geq 1 }$ and a continuous map $\Psi : \Pp_2(\RR^d) \to \RR$,   such that
\begin{align}
    \label{Eq:optimal-policy-convergence}
    \lim_{k \to \infty} \sup_{ x^{n_k} \in B^{n_k} } 
        \Abs{ \Psi^{n_k}( x^{n_k} ) -\Psi( \overline{\mu}_{ x^{n_k} } ) }
    =
    0,
\end{align}
where $\overline{\mu}_{ x^{n} }$ denotes the empirical measure on the vector $x^n = (x^{1,n},\ldots,x^{n,n})$, see \cite[Thm. 2.1, 5.10]{cardaliaguet_notes_2013}.

Thus, for $n$ sufficiently large, the dynamics of $X^n$ in \eqref{Eq:n-player-dynamics} can be approximated (heuristically)  by
\begin{align}
    \label{Eq:n-player-MF-dynamics}
    X^{i,n}_0
    &+
    \int_0^t b \Parentheses{ \overline{\mu}_{ X^{-i,n}_{s-} },  X^{i,n}_{s-}, a^{i,n}_{s-} } \diff s
    +
    \int_0^t \sigma \Parentheses{ \overline{\mu}_{ X^{-i,n}_{s-} }, X^{i,n}_{s-}, a^{i,n}_{s-} } 
        \diff W_s^i 
    \\ \nonumber
    &+
    \int_{[0,t] \times \mathbf{R} } 
        \gamma \Parentheses{ r, \overline{\mu}_{ X^{-i,n}_{s-} }, X^{i,n}_{s-}, a^{i,n}_{s-} }
        \tilde{\eta}^{\lambda,n} (\diff s, \diff r),
\end{align}
where $\tilde{\eta}^{\lambda,n}$ is the compensated measure for the lifting $\eta^{\lambda,n}$ given by the empirical noise $( N^{\lambda,n}$, $\xi^{\lambda,n} )$ with kernel $K^{\lambda,n}_t(\diff r) = \lambda( r, \overline{\mu}_{X^n_{t-}} )Q(\diff r)$ (see Remark \ref{Remark:Form-of-lambda}), such that $b,\sigma,\gamma,\lambda$ are the limits of $b^n,\sigma^n,\gamma^n,\lambda^n$ in the sense of \eqref{Eq:optimal-policy-convergence}. Similarly the payoff $J^n( X^{n}, \alpha^{i,n} )$ from \eqref{Eq:n-player-payoff} can be approximated by $J( \overline{\mu}_{X^{-i,n}}, X^{i,n}, \alpha^{i,n} )$, where
\begin{align}
    \label{Eq:n-player-MF-payoff}
    J( \mu, X, \alpha )
    :=
    \Esp{ \int_0^T f( \mu_{t-}, X_{t-}, \alpha_{t-} ) \diff t + g( \mu_{T-}, X_{T-} ) },
\end{align}
with $f$ and $g$ being the limits of $f^n$ and $g^n$. Lastly, if the $n$-player symmetric game \eqref{Eq:n-player-dynamics}-\eqref{Eq:n-player-payoff} admits a Nash equilibrium with optimal policies $\widehat{\mathbf{a}}^{1,n},\ldots,\widehat{\mathbf{a}}^{n,n}$, they could also be approximated by a continuous function 
\begin{align}
    \label{Eq:n-player-MF-optimal-policy}
    \widehat{\mathbf{a}}: [0,T] \times \Pp_2(\RR^d) \times \RR^d \to A,
\end{align}
defined as in \eqref{Eq:optimal-policy-convergence}. Whenever this type of coupling occurs, i.e. when the interactions between players is only due to the empirical measure or the population, we say that the game has \textit{Mean-Field interactions}.

With these considerations in mind, suppose that the $n$-player symmetric game \eqref{Eq:n-player-dynamics}-\eqref{Eq:n-player-payoff} can be written directly as a game with Mean-Field interactions:
\begin{align}
    \label{Eq:Empirical-Payoff}
    &J^{\overline{\mu}}( X^{i,n}, \alpha^{i,n} ) = J(\overline{\mu}, X^{i,n}, \alpha^{i,n} ),
    &
    &\forall i =1, \ldots, n,
\end{align}
where $\overline{\mu} = \overline{\mu}_{X^n}$ and $J$ is as in \eqref{Eq:Empirical-Payoff}, subject to the dynamics
\begin{align}
    \label{Eq:Empirical-Dynamics}
    X^{i,n}_t
    =&
    X^{i,n}_0
    +
    \int_0^t b^{\overline{\mu}} \Parentheses{ X^{i,n}_{s-}, \alpha^{i,n}_{s-} } \diff s
    +
    \int_0^t \sigma^{\overline{\mu}} \Parentheses{ X^{i,n}_{s-}, \alpha^{i,n}_{s-} } \diff W_s^i 
    \\ \nonumber
    &+
    \int_{[0,t] \times \mathbf{R} } 
        \gamma^{\overline{\mu}} \Parentheses{ r, X^{i,n}_{s-}, \alpha^{i,n}_{s-} }
    \tilde{\eta}^{\lambda,n} (\diff s, \diff r),
\end{align}
defined on the empirical set-up $\big( X_0, \overline{\mu}, K^{\lambda,n}, (N^{\lambda,n},\xi^{\lambda,n}), (W^1,\ldots,W^n) \big)$.

Suppose that the game \eqref{Eq:Empirical-Payoff}-\eqref{Eq:Empirical-Dynamics} admits a Nash equilibrium $(\widehat{X}^n,\widehat{\alpha}^n)$ for an optimal policy $\widehat{\mathbf{a}}$. Then, from \cite[Thm. 9]{hernandez-hernandez_conditional_2023},  under a set of suitable assumptions, we conclude that for any finite set of indexes $I \subset \NN$ the optimally controlled system $\widehat{X}^{I,n}$ converges to a collection of pathwise-exchageable processes $\widehat{X}^{I}$ such that the pair $(\widehat{\mu}^i, \widehat{X}^i) := ( \Ll^1( \widehat{X}^{i} ), \widehat{X}^{i} )$ satisfies the conditional McKean-Vlasov SDE
\begin{align}
    \nonumber
    \widehat{X}^i_t
    =&
    \widehat{X}^{i}_0
    +
    \int_0^t b \Parentheses{ 
        \widehat{\mu}^{i}_{s-}, \widehat{X}^{i}_{s-}, 
        \widehat{\mathbf{a}}_{s-} \big( \widehat{\mu}^{i}_{s-}, \widehat{X}^{i}_{s-} \big)
    } \diff s
    +
    \int_0^t \sigma \Parentheses{ 
        \widehat{\mu}^{i}_{s-}, \widehat{X}^{i}_{s-}, 
        \widehat{\mathbf{a}}_{s-} \big( \widehat{\mu}^{i}_{s-}, \widehat{X}^{i}_{s-} \big)
    } \diff W_s^i 
    \\ \label{Eq:(Example-III)Underline_X}
    &+
    \int_{[0,t] \times \mathbf{R} } \gamma \Parentheses{ 
        r, \widehat{\mu}^{i}_{s-}, \widehat{X}^{i}_{s-}, 
        \widehat{\mathbf{a}}_{s-} \big( \widehat{\mu}^{i}_{s-}, \widehat{X}^{i}_{s-} \big)
    } \tilde{\eta}^{\lambda^i}(\diff s, \diff r),
\end{align}
for each $i \in I$, where $\Ll^1(X)$ denotes the conditional law of $X$ given $(X_0,\mu,\lambda,(N^\lambda,\xi^\lambda))$, $\eta^{\lambda^i}$ is the lifting of the admissible common noise process $(N^{\lambda^i},\xi^{\lambda^i})$ with intensity kernel $K^{\lambda^i}_t(\diff r) = \lambda(r,\widehat{\mu}^i) Q(\diff r)$, and $\tilde{\eta}^{\lambda}$ is the compensated martingale of $\eta^\lambda$. Hence, since for each $i \in I$
\begin{align*}
    \Parentheses{ 
        \Omega, \Ff, \FF, \PP, 
        \widehat{X}^i_0, \widehat{\mu}^i, \lambda^i, (N^{\lambda^i}, \xi^{\lambda^i}), W^i 
    }
\end{align*}
is an admissible set-up, in the sense of Definition \ref{Def:Admissible-Set-Up}, a natural question  is the following: When  each individual expected cost at equilibrium also converges to
\begin{align*}
    J \big( 
        \widehat{\mu}^i, \widehat{X}^i, \widehat{\mathbf{a}}_\cdot( \widehat{\mu}^i, \widehat{X}^i ) 
    \big)
    =
    J \big(
        \widehat{\mu}^1, \widehat{X}^1, \widehat{\mathbf{a}}_\cdot( \widehat{\mu}^1, \widehat{X}^1 ) 
    \big) ?
\end{align*}
In other words, \textit{if the payoff associated with an arbitrary (representative) player can  be interpreted, in the limit, as the payoff of an infinite player game}

This intuition  leads to the next formulation of  the  \emph{Mean-Field Game with common Poissonian noise}, which is formulated in a precise way below.
\begin{definition}
\label{Def:MFE}
The \emph{Mean-Field Game with common Poissonian noise} is defined as the problem consisting on finding an $\FF$-adapted, \cadlag\ process $(\widehat{\mu}, \widehat{X}) = (\widehat{\mu}, \widehat{X}^{\widehat{\mu}})$, defined on some filtered probability space $(\Omega,\Ff,\FF,\PP)$, such that 
\begin{align}
    \label{Eq:Fixed-Point-Problem}
    \widehat{\mu}
    \mbox{ is a fixed point of the mapping }
    \mu \mapsto \Ll^1(\widehat{X}^\mu),
\end{align}
where $\widehat{X}^\mu$ is defined as the optimal solution to the \emph{control problem in random environment $\mu$}, defined on an admissible set-up $( X_0, \mu, \lambda, (N^\lambda,\xi), W)$, given by
\begin{align}
    \label{Eq:MFG-Payoff}
    \underset{\alpha \in \Aa}{\mathrm{minimize}}\  J^\mu(X,\alpha),\;
    \mbox{ with}\;\;
    J^\mu(X,\alpha) := \Esp{ \int_0^T f^\mu ( X^\mu_{t-}, \alpha_{t-} ) \diff t + g^\mu( X^\mu_{T-} ) },
\end{align}
subject to the dynamics
\begin{align}
    \label{Eq:MFG-System}
    X^\mu_t
    =&
    X_0
    +
    \int_0^t b^\mu \Parentheses{ X^\mu_{s-}, \alpha_{s-} } \diff s
    +
    \int_0^t \sigma^\mu \Parentheses{ X^\mu_{s-}, \alpha_{s-} } \diff W_s
    \\ \nonumber
    &+
    \int_{[0,t] \times \mathbf{R}} \gamma^\mu \Parentheses{ r, X^\mu_{s-}, \alpha_{s-} }
        \tilde{\eta}^\lambda(\diff s, \diff r).
\end{align}
We refer to $(\widehat{\mu}, \widehat{X})$ as the \emph{Mean-Field Equilibrium (MFE) for the game}. If  $(\widehat{\mu}, \widehat{X})$ exists and  is defined on the same admissible set-up
\begin{align*}
    \Parentheses{ \Omega, \Ff, \FF, \PP, X_0, \Ll^1(\widehat{X}), \lambda, (N^\lambda, \xi), W },
\end{align*}
we say that $(\widehat{\mu}, \widehat{X})$ is a \emph{strong MFE}.
\end{definition}

\begin{remark}
In the following, we refer to the controlled pair $(\widehat{X},\widehat{\alpha})$ as the \textit{optimal pair of the MFG} whenever $(\widehat{\mu},\widehat{X})$ is a MFE for \eqref{Eq:MFG-Payoff}, and $\widehat{X}$ solves the $\widehat{\alpha}$-controlled Conditional McKean-Vlasov Equation \eqref{Eq:MFG-System} with $\mu = \widehat{\mu}$. To ease the notation, we write 
\begin{align}
    \label{Eq:Hat-notation}
    \widehat{b}_t(\omega)
    &:=
    b( \widehat{\mu}_{t-}(\omega), \widehat{X}_{t-}(\omega), \widehat{\alpha}_{t-}(\omega)),
    &
    &\diff \PP \otimes \diff t-\mbox{a.e.}
\end{align}
(and similarly for $\sigma,\gamma,\lambda,\ldots$) when evaluating $\widehat{\alpha}$-controlled systems under the equilibrium $(\widehat{\mu},\widehat{X})$, with the differences being represented by 
\begin{align*}
    &\delta \widehat{b}_t^\nu( x, a )(\omega) 
    = 
    {b}_t^\nu( x, a )(\omega) - \widehat{b}_t(\omega),
    &
    &\diff \PP \otimes \diff t-\mbox{a.e.},
\end{align*}
for all $(\nu, x, a) \in \Pp(\RR^d) \times \RR^d \times A$; see Section \ref{Sec:Basic-Notation} above.

\end{remark}

Now that the game and its equilibrium have been defined, our goal is to explore a way to solve it or, at least,  find a set of conditions that guarantee the existence of a solution. However, recalling the results on propagation of chaos shown in  \cite{hernandez-hernandez_coupled_2023}, we notice that McKean-Vlasov equations necessarily solve the fixed point problem described in \eqref{Eq:Fixed-Point-Problem}. Indeed, if $(\widehat{\mu}, \widehat{X})$ is a MFE with the corresponding control $\widehat{\alpha}$, then the dynamics of $\widehat{X}$ satisfy the McKean-Vlasov equation
\begin{align*}
    \widehat{X}_t
    =&
    X_0
    +
    \int_0^t b \Parentheses{ \Ll^1(\widehat{X})_{s-}, \widehat{X}_{s-}, \widehat{\alpha}_{s-} } \diff s
    +
    \int_0^t \sigma \Parentheses{ \Ll^1(\widehat{X})_{s-}, \widehat{X}_{s-}, \widehat{\alpha}_{s-} } 
        \diff W_s
    \\ \nonumber
    &+
    \int_{[0,t] \times \mathbf{R}} 
        \gamma \Parentheses{ r, \Ll^1(\widehat{X})_{s-}, \widehat{X}_{s-}, \widehat{\alpha}_{s-} }
        \tilde{\eta}^\lambda(\diff s, \diff r).
\end{align*}
With this in mind, we shall focus next in 
solving the optimization problem described in \eqref{Eq:MFG-Payoff}-\eqref{Eq:MFG-System}.  
This can be seen as an optimal control problem in random environment, according with the description given in Definition \ref{Def:MFE}.

\section{Pontryagin's maximum principle}
\label{Sec:PMP}

Our approach  is based on the so-called \textit{Pontryagin's Maximum Principle}; namely, we use variational arguments to develop a set of necessary conditions that an optimal trajectory must fulfil, as well an adjoint set of differential equations which capture the variation of the payoff functional $J^\mu(X,\alpha)$ when the control $\alpha$  is slightly perturbed. It should be emphasized that throughout this section the random environment $\mu$ is arbitrary but fixed. Technically, our method is based on the works of Tang \& Li \cite{tang_necessary_1994} and  Framstad, {\O ksendal} \& Sulem \cite{framstad_sufficient_2004}, which considered systems driven by jump-diffusion, but on a considerably simpler setting, i.e. with no random environment.

\subsection{Hamiltonian system: Definitions and Assumptions}

The method presented in this paper to verify the PMP is based in \cite{tang_necessary_1994},  and systematically use variational arguments in order to obtain the Hamiltonian system corresponding to the game problem  in random environment. This means that, given an optimal pair $(\widehat{X},\widehat{\alpha})$ for the optimization problem in \eqref{Eq:MFG-Payoff}-\eqref{Eq:MFG-System}, our goal is to derive the set of adjoint backward SDEs which captures the first and second order variations of the optimally controlled system $\widehat{X}$, given a small perturbation of the optimal control $\widehat{\alpha}$. The Hamiltonian, defined below, is a crucial tool for this  purpose, in order to describe the variations in the total cost  $J^\mu$ when a representative player  deviates from the optimally controlled trajectory $\widehat{X}$.

The \emph{$\lambda$-admissible Hamiltonian} (or simply \emph{Hamiltonian}) is defined as a mapping
\begin{align*}
    H
    \ : \ 
    \Pp_2(\RR^d) \times \RR^d \times \RR^d \times \RR^{d \times k} \times \mathfrak{H}^\lambda \times A 
    \to
    \RR
\end{align*}
given by
\begin{align}
    \label{Eq:H}
    H(\nu,x,y,z,u,a)
    :=&
    f(\nu,x,a)
    +
    \Parentheses{ y \cdot b(\nu,x,a) }
    +
    \Parentheses{ z \cdot \sigma(\nu,x,a) }_\mathrm{Fr}
    \\ \nonumber
    &+
    \int_\mathbf{R} \Tr{ u(r) \diag{\lambda(r,\nu)Q(\diff r)} \gamma(r,\nu,x,a)^\top },
\end{align}
on the event $\Braces{ U_t = u }$ for $U \in \mathfrak{H}^\lambda$.

The following standing hypotheses are assumed for the rest of the paper. We start with a couple of technical conditions on the setting: the first two determine the admissiblity of a common Poissonian noise and its form, while the third one is related to the existence of the required \textit{a priori} estimates. As for  the coefficients, we introduce  conditions which ensure that the Hamiltonian system (described below in terms of the adjoint FBSDE) is well-posed \cite{hernandez-hernandez_coupled_2023}, \cite{hernandez-hernandez_conditional_2023}. Here, recall that we are using the $\delta x$ notation to denote the difference $x - x'$. Finally, we conclude with two technical conditions to ensure Pontryagin's (maximal) condition.
\begin{assumptions}
\label{Assumptions-C}
\hfill
\par
\noindent
\textbf{(i) On the probability space.} 
 
\begin{enumerate}
    \renewcommand{\labelenumi}{(C\theenumi)}
    

        
        

    \item
        The tuple $(\Omega, \Ff, \PP, X_0, \mu, \lambda, (N^\lambda,\xi),W)$ is a given admissible set-up such that $(\Omega,\Ff)$ is a Polish space endowed with its Borel $\sigma$-algebra.
        
    \item
        The intensity candidate has the form $\lambda_t(\cdot,r) = \lambda(r, \mu_{t-}(\cdot))$ for some measurable deterministic function $\lambda:\mathbf{R} \times \Pp_2(\RR^d) \to \RR_+^l$. Furthermore, $\lambda$ is a Lipschitz function on $\Pp_2(\RR^d)$ with respect to the $2$-Wasserstein distance.

    \item 
        The random environment $\mu$ satisfies the uniform bound
        \begin{align*}
            \Norm{ \mu }_\Ww 
            := 
            \sup_{[0,T]} \Parentheses{ \Esp{ W_2\Parentheses{ \mu_t, \delta_{ \{ 0 \} } }^8 } }^\frac{1}{8}
            < 
            \infty.
        \end{align*}

    \setcounter{assumption-counter}{\value{enumi}}
\end{enumerate}


\noindent
\textbf{(ii) On the coefficients.}  
\begin{enumerate}
    \renewcommand{\labelenumi}{(C\theenumi)}
    \setcounter{enumi}{\value{assumption-counter}}
    
    \item
        \emph{Second differentiability on the state.}
        \begin{itemize}
            \item 
                $(b,\sigma,f,g)(\nu, \cdot, a) \in \mathrm{C}_b^2(\RR^d)$ for all $( \nu, a ) \in \Pp_2(\RR^d) \times A$;

            \item 
                $\gamma( r, \nu, \cdot, a) \in \mathrm{C}_b^2(\RR^d)$ for all $( r, \nu, a ) \in \mathbf{R} \times \Pp_2(\RR^d) \times A$;

            \item 
                $g(\nu, \cdot) \in \mathrm{C}_b^2(\RR^d)$ for all $\nu \in \Pp_2(\RR^d)$.
                
        \end{itemize}

    \item
        \emph{Uniform continuity on the control.}
        \begin{itemize}
            \item 
                $( \Psi, \Diff \Psi, \Diff^2 \Psi )( \nu, x, \cdot) \in \mathrm{UC}(A)$ for all $(\nu, x) \in \Pp_2(\RR^d) \times \RR^d$, where $\Psi = b, f, \sigma$;

            \item 
                $( \gamma, \Diff \gamma, \Diff^2 \gamma )( r, \nu, x, \cdot ) \in \mathrm{UC}(A)$ for all $( r, \nu, x ) \in \mathbf{R} \times \Pp_2(\RR^d) \times \RR^d$.
                
        \end{itemize}

    \item
        \emph{Lipschitz condition on the measure and the adjoint states.} \newline
        There exists a positive constant $K > 0$, such that for any fixed $a \in A$, the following inequality holds for $\diff \PP \otimes \diff t$-almost all $(\omega,t) \in \Omega \times [0,T]$:
		\begin{align*}
            &\Abs{ \delta b( \nu, x \vert a ) }^2
            +
            \Norm{ \delta \sigma( \nu, x \vert a ) }^2_{\mathrm{Fr}}
            \\
            &\hspace{3em}+
            \int_{\mathbf{R}} \Tr{ 
                \delta \gamma(\nu,x \vert r,a) 
                \diag{ \lambda(r,\nu) Q(\diff r) } 
                \delta \gamma(\nu,x \vert r, a)^\top 
            }
            \\
            &\hspace{3em}+
            \Abs{ \delta \Diff H (\nu, x, y, z, u \vert a ) }^2
            + 
            \Abs{ \delta \Diff g(\nu,x) }^2
            \\
            &\leq
            K \bigg\lbrace 
                W_2(\nu,\nu')^2 
                +
                \Abs{ \delta x }^2 
                + 
                \Abs{ \delta y }^2 
                + 
                \Norm{ \delta z }_{ \mathrm{Fr} }^2 
                \\
                &\hspace{3em}+ 
                \int_{\mathbf{R}} \Tr{ 
                    \delta u(r) \diag{ \lambda(r,\nu) Q(\diff r) } \delta u( r )^\top 
                }
            \bigg\rbrace,
		\end{align*}
        with $H$ as  in \eqref{Eq:H}, for each 
        \begin{align*}
            ( \nu, x, y, z, u), ( \nu', x', y', z', u' )  \in \Pp_2(\RR^d) \times \RR^d \times \RR^d \times \RR^{ d \times m } \times  \mathfrak{H}^\lambda
        \end{align*} 
        on the event $\Braces{ u_t(\cdot) - u'_t(\cdot) = \delta u(\cdot) }$.
        
    \item
        \emph{Monotonicity condition on the adjoint states.} \newline
        For any fixed $\nu \in \Pp_2(\RR^d)$ and $a \in A$, there exist nonnegative constants $\beta_1$, $\beta_2$, $\beta_3$ with $\beta_1 + \beta_2 > 0$ and $\beta_2 + \beta_3 > 0$ such that the following inequalities hold $\diff \PP \otimes \diff t$-a.e.:
        \begin{align*}
            &\Big( \delta \Diff H ( x, y, z, u \vert \nu, a ) \cdot \delta x \Big) 
            +
            \Big( \delta b ( x \vert \nu, a ) \cdot \delta y \Big) 
            +
            \Big( \delta \sigma ( x \vert \nu, a ) \cdot \delta z \Big)_\mathrm{Fr}
            \\
            &\hspace{10em}+
            \int_{\mathbf{R}} \Tr{ 
                \delta \gamma( x \vert r, \nu, a ) 
                \diag{ \lambda(r,\nu) Q(\diff r) } 
                \delta u(r)^\top 
            }
            \\
            &\leq
            - 
            \beta_1 \Abs{ \delta x }^2
            -
            \beta_2 \bigg\lbrace 
                \Abs{ \delta y }^2 + \Norm{ \delta z }^2_{\mathrm{Fr}} 
                + 
                \int_{\mathbf{R}} \Tr{ \delta u(r) \diag{ \lambda(r,\nu) Q(\diff r) } \delta u( r )^\top }
            \bigg\rbrace,
        \end{align*}
        and
        \begin{align*}
            &\Parentheses{ \delta g( x \vert \nu ) \cdot \delta x }
            \geq 
            \beta_3 \Abs{ \delta x }^2,
        \end{align*}
        for all $ (x,y,z,u), (x',y',z',u') \in \RR^d \times \RR^n \times \RR^{ d \times k } \times \mathfrak{H}^\lambda$ on the event $\{ u_t(\cdot) - u'_t(\cdot) = \delta u(\cdot) \}$.

    \item 
        \emph{Moment condition on the marks.}\newline
        There exists a positive constant $\gamma^*$ such that 
        \begin{align*}
            &\max_{i \in \{1,\ldots,d\} } 
            \int_\mathbf{R} \Abs{ \gamma^{(i,\cdot)}(r,\nu,x,a)}^4 Q^{(i)}( \diff r ) 
            \leq \gamma^* < \infty,
        \end{align*}
        for all $(\nu,x,a) \in \Pp_2(\RR^d) \times \RR^d \times A$.

    \setcounter{assumption-counter}{\value{enumi}}
\end{enumerate}


\noindent
\textbf{(iii) On the controls:}  
\begin{enumerate}
    \renewcommand{\labelenumi}{(C\theenumi)}
    \setcounter{enumi}{\value{assumption-counter}}
    \item 
        The set of admissible actions $A$ is given by a non-empty, convex and compact subset of a finite-dimensional Euclidean space, endowed with its usual topology.

    \item
        The set of admissible strategies is given by 
        \begin{align*}
            \mathcal{A}
            :=
            \Big\lbrace
                \alpha : \Omega \times [0,T] \to A
                \ \Big\vert \ 
                &\alpha\mbox{ is \cadlag\ and }\FF-\mbox{adapted, with}
                \\
                &\hspace{3em}\Norm{ \alpha }_\Aa
                :=
                \sup_{[0,T]} \Parentheses{ \Esp{ \Abs{ \alpha_t }^8 } }^\frac{1}{8}
                <
                \infty
            \Big\rbrace.
        \end{align*}

\end{enumerate}

\end{assumptions}

In order to derive the first and second order variations of the optimally controlled trajectory, we  need to analyze  the integral term appearing in the usual Taylor's expansion. That is, if $\Psi \in \mathrm{C}^p(V)$, with $p \in \NN$ and $V \subset \RR^d$ an open and convex set, then
\begin{align}
    \label{Eq:Taylor-expansion}
    &\Psi(x) - \sum_{k=0}^{p-1}\frac{1}{k!} \Diff^{k}\Psi(x_0;x-x_0)
    \\ \nonumber
    &{\hskip0.15\linewidth}
    =
    \frac{1}{(p-1)!} \int_0^1 
        \Diff^{p}\Psi\Parentheses{ x_0 + \tau(x-x_0); x - x_0 } \diff \tau,
\end{align}
for all $x,x_0$ in the interior of $V$, where $\Diff^{k}\Psi(x_0;x-x_0)$ denotes the $k$-th differential of $\Psi$ around $x_0$. 
\begin{remark}
Given a random environment $\mu$, hereafter we write $H^\mu = \{ H_t^\mu, 0 \leq t \leq T \}$ to denote the random process
\begin{align*}
    &H_t^\mu(x,y,z,u,a)(\omega) := H(\mu_t(\omega), x, y, z, u, a ), 
    &
    &\diff \PP \otimes \diff t-\mbox{a.e.}
\end{align*}
for  $(x,y,z,u,a) \in \RR^d \times \RR^d \times \RR^{d \times k} \times \mathfrak{H}^\lambda \times A$, with $\delta H_t^\mu$ and $\widehat{H}_t^\mu$ following the notation stated in  Section \ref{Sec:Basic-Notation} and equation \eqref{Eq:Hat-notation}, respectively. We also use $(v \cdot M)$ to refer to the vector-matrix multiplication $v^\top M$ whenever the dimensions coincide and the product is well-defined.
\end{remark}

Observe that 
the boundedness $\diff \PP \otimes \diff t$-a.a. on the derivatives and the conditions imposed on the intensity process imply that  the Jacobian of $H$ with respect to $x$ is given by
\begin{align}
    \label{Eq:Diff_H}
    \Diff H_t^\mu( x, y, z, u, a )
    =&
    \Diff f_t^\mu( x, a )
    +
    \Parentheses{y \cdot \Diff b_t^\mu( x, a )}
    +
    \sum_{i=1}^k \Parentheses{ z^{(\cdot,i)} \cdot {\Diff \sigma^{(\cdot,i)}}_t^\mu( x, a )}
    \\ \nonumber
    &+
    \sum_{j=1}^l \Parentheses{ 
        u_t^{(\cdot,j)} \cdot {\Diff \gamma^{(\cdot,j)}}_t^\mu(\cdot,x,a) 
    }_{ \lambda^{(j)}_t }.
\end{align}
Notice that the last term in the above display involves an integral with respect to the random measure, and equalities  should be understood   component-wise for $\diff\PP \otimes \diff t$-a.a. $(\omega,t) \in \Omega \times [0,T]$. Similarly, the Hessian of $H$ is given by
\begin{align}
    \label{Eq:Diff^2_H}
    \Diff^2{H}_t^\mu( x, y, z, u, a )
    =&
    \Diff^2{f}_t^\mu( x, a )
    +
    \Big( \Parentheses{ y \otimes \Id_{d} } \cdot \Diff^2{b}_t^\mu( x, a ) \Big)
    \\[3pt] \nonumber
    &+
    \sum_{i=1}^k \Parentheses{
        \Parentheses{ z^{(\cdot,i)} \otimes \Id_d } 
        \cdot 
        {\Diff^2\sigma^{(\cdot,i)}}_t^\mu( x, a )
    }
    \\ \nonumber
    &+
    \sum_{j=1}^l \Parentheses{
        \big( u_t^{(\cdot,j)} \otimes \Id_d \big)
        \cdot
        {\Diff^2\gamma^{(\cdot,j)}}_t^\mu(\cdot,x,a)
    }_{ \lambda^{(j)}_t },
\end{align}
where $\Id_d$ represents the $d$-dimensional identity matrix.

\begin{remark}
Some  computations presented below require  algebraic properties of the Kroenecker product, which are recalled next
\begin{align*}
    (\Id_d \otimes v) v_0
    =&
    (\Id_d \otimes v) \mathrm{vec}(v_0^\top)
    =
    \mathrm{vec}( v v_0^\top \Id_d )
    \\
    =&
    \mathrm{vec}( \Id_d v v_0^\top )
    =
    (v_0 \otimes \Id_d) \mathrm{vec}(v)
    =
    (v_0 \otimes \Id_d) v
\end{align*}
for  $v,v_0 \in \RR^d$. They imply, in particular, that  for any $M \in \RR^{d \times d}$
\begin{align*}
    v_0^\top \Big( \Parentheses{ \Id_{d} \otimes v } \cdot M v \Big)
    =&
    \Big( \Parentheses{ \Id_{d} \otimes v }v_0 \cdot M v \Big)
    \\
    =&
    \Big( \Parentheses{ v_0 \otimes \Id_{d} }v \cdot M v \Big)
    =
    v^\top \Big( \Parentheses{ v_0 \otimes \Id_{d} } \cdot M \Big) v.
\end{align*}
\end{remark}

Let us assume now that $(\widehat{X},\widehat{\alpha})$ is an optimal pair for the  optimization problem  \eqref{Eq:MFG-Payoff}-\eqref{Eq:MFG-System}. Then, the \textit{first order adjoint process} is defined as the solution  
\begin{align}
    \label{Eq:YZUM}
    (Y,Z,U,M) 
    \in
    \Ss^2(\RR^d) \times \HH^{2}(W) \times \HH^{2}(\tilde{\eta}^\lambda) \times \Hh^{2,\perp}
\end{align}
to the following BSDE:  
\begin{align}
    \label{Eq:Y}
    Y_t
    :=&
    \Diff\widehat{g}_T^{\mu\top}
    +
    \int_t^T \Diff \widehat{ H }_s^{\mu}(Y,Z,U)^\top \diff s
    -
    \int_t^T Z_s \diff W_s
    \\ \nonumber
    &-
    \int_{ (t,T] \times \mathbf{R} } U_{s}(r) \tilde{\eta}^\lambda (\diff s, \diff r)
    -
    (M_T - M_t).
\end{align}
Similarly, the \textit{second order adjoint process} 
\begin{align}
    \label{Eq:Tildes-YZUM}
    \Parentheses{\widetilde{Y},\widetilde{Z},\widetilde{U},\widetilde{M}}
    \in
    \Ss^2(\RR^{d \times d}) \times \HH^{2}(W)^k \times \HH^{2}(\tilde{\eta}^\lambda)^l \times \Hh^{2,\perp},
\end{align}
with $\widetilde{Z} =(\widetilde{Z}^{(1)},\ldots,\widetilde{Z}^{(k)})$ and $\widetilde{U} =(\widetilde{U}^{(1)},\ldots,\widetilde{U}^{(l)})$, is given as the solution to the BSDE
\begin{align}
    \label{Eq:tilde-Y}
    \widetilde{Y}_t
    :=&
    \Diff^2 \widehat{ g }_T^\mu
    +
    \int_t^T B^{\star\star}_s \diff s
    -
    \sum_{i=1}^k \int_t^T \Big\lbrace 
        \widetilde{Z}^{(i)}_s \diff W^{(i)}_s - \Sigma_s^{\star\star(i)} \diff s 
    \Big\rbrace
    \\ \nonumber
    &-
    \sum_{j=1}^l \int_{ (t,T] \times \mathbf{R} } \Big\lbrace
        \widetilde{U}^{(j)}_s(r) \tilde{\eta}^{\lambda(j)}( \diff s, \diff r )
        -
        \Gamma_{s-}^{\star\star(j)}(r) \big(\diff s \otimes K^{\lambda(j)}_s (\diff r) \big)
    \Big\rbrace
    \\[8pt] \nonumber
    &-
    (\widetilde{M}_T - \widetilde{M}_t),
\end{align}
where $\widetilde{Y}$, $\widetilde{Z}^{(i)}$, $\widetilde{U}^{(j)}$, $\widetilde{M}$ are matrix-valued processes of dimension $d \times d$ for each $i$ and  $j$, while $B^{\star\star}, \Sigma^{\star\star(i)}$ and $\Gamma^{\star\star(j)}$ are defined as
\begin{align}
    \label{Eq:B^starstar}
    &B^{\star\star}_t
    := 
    \Diff^2 \widehat{H}_t^\mu( Y, Z, U )
    +
    \Diff \widehat{b}_t^{\mu\top} \widetilde{Y}_t
    +
    \widetilde{Y}_t \Diff \widehat{b}_t^\mu,
    \\ \label{Eq:Sigma^starstar}
    &\Sigma^{\star\star(i)}_t
    :=
    \Diff \widehat{ \sigma }_t^{ \mu (\cdot,i)\top } \widetilde{Y}_t
        \Diff \widehat{ \sigma }_t^{ \mu (\cdot,i) }
    +
    \Diff \widehat{ \sigma }_t^{ \mu (\cdot,i) \top } \widetilde{Z}^{(i)}_t 
    +
    \widetilde{Z}^{(i)}_t \Diff \widehat{ \sigma }_t^{ \mu (\cdot,i) }
    \\ \label{Eq:Gamma^starstar}
    &\Gamma_t^{\star\star(j)}(r)
    :=
    \Diff \widehat{ \gamma }^{\mu (\cdot,j)}_{t}(r)^\top \widetilde{Y}_{t} 
        \Diff \widehat{ \gamma }^{\mu (\cdot,j)}_{t}(r)
    +
    \Diff \widehat{ \gamma }^{\mu (\cdot,j)}_{t}(r)^\top \widetilde{U}^{(j)}_{t}(r)
    \\ \nonumber
    &{\hskip0.15\linewidth}+
    \widetilde{U}^{(j)}_{t}(r) \Diff \widehat{ \gamma }^{\mu (\cdot,j)}_{t}(r),
\end{align}
with the equalities understood component-wise. Then, the \textit{Hamiltonian system for the problem} corresponds to the optimally controlled state trajectory, in conjunction with the \textit{first} and \textit{second order adjoint processes}; in other words, it is defined as the process
\begin{align}
    \label{Eq:Hats-XYZUM}
    \Parentheses{\widehat{X},\widehat{Y},\widehat{Z},\widehat{U},\widehat{M}}
    :=
    \Parentheses{
        \widehat{X}, 
        (Y, \widetilde{Y}),
        (Z, \widetilde{Z}), 
        (U, \widetilde{U}), 
        (M, \widetilde{M})
    },
\end{align}
such that 
\begin{itemize}
    \item 
        $(\widehat{X},\widehat{\alpha}) \in \Ss^2(\RR^{d}) \times \Aa$ is an optimal pair for \eqref{Eq:MFG-Payoff}-\eqref{Eq:MFG-System} under the random environment $\mu$; 
    
    \item 
        $(Y,Z,U,M)$, defined as in \eqref{Eq:YZUM}, solves the \textit{first order variational equation} \eqref{Eq:Y};

    \item 
        $(\widetilde{Y},\widetilde{Z},\widetilde{U},\widetilde{M})$, defined as in \eqref{Eq:Tildes-YZUM}, solves the \textit{second order variational equation} \eqref{Eq:tilde-Y}.
\end{itemize}

\subsection{Necessary conditions for optimality}

Now we proceed to find  necessary conditions for the pair $(\widehat{X},\widehat{\alpha})$, which  must be satisfied in order to be  optimal. These results require some preliminary technical estimations presented next.

\subsubsection{A priori estimates}

\begin{definition}
\label{Def:spike-variation}
Given $\epsilon > 0$, let $I^\epsilon$ be a sub-interval of $[0,T]$ with Lebesgue measure $\Abs{I^\epsilon} = \epsilon$. The \emph{$\epsilon$-variation} (or \emph{spike variation}) \emph{of $\widehat{\alpha}$ with respect to $\alpha$}, denoted by $\alpha^\epsilon \equiv \alpha^{I^\epsilon} \in \Aa$, is then defined as the (admissible) policy
\begin{align}
    \label{Eq:alpha^epsilon}
    \alpha^\epsilon_t(\omega)
    :=
    \begin{cases}
    \alpha_t(\omega),	& t \in I^\epsilon,
    \\
    \widehat{\alpha}_t(\omega),	& t \in [0,T] \setminus I^\epsilon,
    \end{cases}
\end{align}
for all $(\omega,t) \in \Omega \times [0,T]$. Consequently, the $\epsilon$-variational pair is the corresponding admissible pair $(X^\epsilon,\alpha^\epsilon)$
\end{definition}

Observe that the difference between the optimally controlled state $\widehat{X}$ and its $\epsilon$-variation $X^\epsilon$ is described by
\begin{align}
    \label{Eq:delta-hat-X^epsilon}
    \delta \widehat{X}^\epsilon_t
    :=&
    \delta \widehat{X}^\epsilon_0
    +
    \int_{ [0,t] \cap I^\epsilon } \delta \widehat{b}_s^\mu( X, \alpha ) \diff s
    +
    \int_{[0,t] \cap I^\epsilon } \delta \widehat{\sigma}_s^\mu( X, \alpha ) \diff W_s
    \\ \nonumber
    &+
    \int_{ ([0,t] \cap I^\epsilon) \times \mathbf{R} } 
        \delta \widehat{\gamma}_s^\mu( X, \alpha \vert r ) 
    \tilde{\eta}^\lambda (\diff s, \diff r),
\end{align}
with $\delta \widehat{X}^\epsilon_0 = X^\epsilon_0 - \widehat{X}_0$. Following the same intuition behind \cite{yong_stochastic_1999} and \cite{tang_necessary_1994}, the first order state variation of $\delta \widehat{X}^\epsilon$ due to $\alpha^\epsilon$, denoted by $\phi$, follows the dynamics
\begin{align}
    \label{Eq:phi}
    \phi_t
    =&
    \int_0^t \Diff \widehat{ b }_s^\mu \phi_s \diff s
    +
    \sum_{i=1}^k \int_0^t \Braces{
        \delta \widehat{\sigma}^{\mu(\cdot,i)}_s (\alpha^\epsilon \vert \widehat{X} )
        +
        \Diff \widehat{\sigma}^{\mu(\cdot,i)}_s \phi_s
    } \diff W^{(i)}_s
    \\ \nonumber
    &+
    \sum_{j=1}^l \int_{ [0,t] \times \mathbf{R} } \Braces{
        \delta \widehat{\gamma}^{\mu(\cdot,j)}_{s}(\alpha^\epsilon \vert r, \widehat{X} )
        +
        \Diff \widehat{ \gamma }^{\mu(\cdot,j)}_s(r) \phi_{s-}
    } \tilde{\eta}^{\lambda(j)}( \diff s, \diff r ),
\end{align}
for all $t \in [0,T]$, where the equality holds $\PP$-a.s. entry-wise, with 
\begin{align}
    \label{Eq:Hat-notation-2}
    &\delta \widehat{\sigma}^{\mu(\cdot,i)}_t (\alpha^\epsilon \vert \widehat{X} ) 
    = 
    \sigma^{\mu(\cdot,i)}_t (  \widehat{X}, \alpha^\epsilon ) - \widehat{\sigma}^{\mu(\cdot,i)}_t
    \\ \nonumber
    &\Big(\mbox{resp. }
    \delta \widehat{\gamma}^{\mu(\cdot,j)}_t (\alpha^\epsilon \vert r, \widehat{X} ) 
    = 
    \gamma^{\mu(\cdot,j)}_t ( r,  \widehat{X}, \alpha^\epsilon ) 
        - \widehat{\gamma}^{\mu(\cdot,j)}_t(r)
    \Big),
\end{align}
in accordance to the notation from \eqref{Eq:Hat-notation}, and
\begin{align*}
    &\Diff \widehat{\sigma}^{\mu(\cdot,i)}_t
    =
    \begin{bmatrix}
        \Diff \widehat{\sigma}^{\mu(1,i)}_t
        \\
        \vdots
        \\
        \Diff \widehat{\sigma}^{\mu(d,i)}_t
    \end{bmatrix},
    &
    &\Diff \widehat{\gamma}^{\mu(\cdot,j)}_t(\cdot)
    =
    \begin{bmatrix}
        \Diff \widehat{\gamma}^{\mu(1,j)}_t(\cdot)
        \\
        \vdots
        \\
        \Diff \widehat{\gamma}^{\mu(d,j)}_t(\cdot)
    \end{bmatrix},
\end{align*}
for each $i = 1, \ldots, k$ and $j = 1, \ldots,l$, respectively.

Analogously, the second order variation $\vartheta$ is described as the solution to the SDE
\begin{align}
    \label{Eq:vartheta}
    \vartheta_t
    =&
    \delta \widehat{X}^\epsilon_0
    +
    \int_0^t \Big\lbrace
		\delta \widehat{b}_s^\mu(\alpha^\epsilon \vert \widehat{X} )
		+
		\Diff \widehat{b}_s^\mu \vartheta_s 
		+
		\frac{1}{2}
            \Parentheses{ \Id_d \otimes \phi_s }^\top \Diff^2 \widehat{ b }_s^\mu \phi_s
    \Big\rbrace \diff s
    \\ \nonumber
    &+
    \sum_{i=1}^k \int_0^t \Big\lbrace
        \Diff \widehat{ \sigma }_s^{\mu (\cdot,i)} \vartheta_s
        +
        \delta \Diff \widehat{ \sigma }_s^{\mu (\cdot,i)}(\alpha^\epsilon \vert \widehat{X}) \phi_s
        +
        \frac{1}{2}
            (\Id_d \otimes \phi_s )^\top \Diff^2 \widehat{ \sigma }_s^{\mu (\cdot,i)} \phi_s
    \Big\rbrace \diff W^{(i)}_s
    \\ \nonumber
    &+
    \sum_{j=1}^l \int_{ [0,t] \times \mathbf{R} } \Big\lbrace
        \Diff \widehat{ \gamma }_{s-}^{\mu (\cdot,j)}(r) \vartheta_{s-}
        +
        \delta \Diff \widehat{ \gamma }_{s-}^{\mu (\cdot,j)}(\alpha^\epsilon \vert r, \widehat{X}) \phi_{s-}
        \\ \nonumber
        &\qquad\qquad\qquad\quad+
        \frac{1}{2}
            (\Id_d \otimes \phi_{s-} )^\top 
            \Diff^2 \widehat{ \gamma }_{s-}^{\mu (\cdot,j)}(r) 
            \phi_{s-}
    \Big\rbrace \tilde{\eta}^{\lambda(j)}( \diff s, \diff r ),
\end{align}
with $\delta \Diff \widehat{\sigma}^{\mu(\cdot,i)}_t (\alpha^\epsilon \vert \widehat{X} )$ and $\delta \Diff \widehat{\gamma}^{\mu(\cdot,j)}_t (\alpha^\epsilon \vert r, \widehat{X} )$ following a similar notation to \eqref{Eq:Hat-notation-2}, and
\begin{align*}
    &\Diff^2 \widehat{\sigma}^{\mu(\cdot,i)}_t
    =
    \begin{bmatrix}
        \Diff^2 \widehat{\sigma}^{\mu(1,i)}_t
        \\
        \vdots
        \\
        \Diff^2 \widehat{\sigma}^{\mu(d,i)}_t
    \end{bmatrix},
    &
    &\Diff^2 \widehat{\gamma}^{\mu(\cdot,j)}_t(\cdot)
    =
    \begin{bmatrix}
        \Diff^2 \widehat{\gamma}^{\mu(1,j)}_t(\cdot)
        \\
        \vdots
        \\
        \Diff^2 \widehat{\gamma}^{\mu(d,j)}_t(\cdot)
    \end{bmatrix},
\end{align*}
for each $i = 1, \ldots, k$ and $j = 1, \ldots,l$, respectively.

\begin{remark}
\hfill
\noindent
\begin{itemize}
    \item 
        Observe that from the regularity imposed on the coefficients, i.e. Assumptions \ref{Assumptions-C}.(ii), the linear SDEs described in \eqref{Eq:phi} and \eqref{Eq:vartheta} have a strong solution, (see the corresponding discussion on uncoupled linear components in \cite{hernandez-hernandez_coupled_2023}).

    \item 
        As in the classical case, notice that when the control does not appear in either the diffusion or the jump terms, we can  discard $\phi$ from them, since it does not contribute to the variation of the control.
\end{itemize}

\end{remark}

The processes $\phi$ and $\vartheta$ above correspond to the stochastic version of the estimates shown in \cite{tang_necessary_1994}. From \cite{hernandez-hernandez_coupled_2023}, we know the coefficients involved require additional measurability conditions in order to guarantee the existence of both processes. Note that these conditions are satisfied once the well-posedeness of the controlled state $X^\alpha$ is guaranteed. However, in order to attain the appropriate estimates for the variation of the cost, first we require to control the randomness inherent to the environment.

\begin{lemma}
\label{Lemma:A-Priori-Estimates}
Let $(\widehat{X},\widehat{\alpha})$ be an optimal pair, and for $\epsilon > 0$ let $(X^\epsilon,\alpha^\epsilon)$ be the corresponding spike variation with respect to $\alpha \in \Aa$. 
Then, for each $k=1,\ldots,4$,
\begin{align}
    \label{Eq:A-Priori-I}
    &\sup_{[0,T]} \Esp{ \Abs{ X_t^\alpha }^{2k} }
    \leq
    O\Parentheses{ (1 + \Norm{\mu}_\Ww + \Norm{\alpha}_\Aa )^{2k} },
    \\
    \label{Eq:A-Priori-II}
    &\sup_{[0,T]} \Esp{ \Abs{ \delta \widehat{X}^\epsilon_t }^{2k} }
    \leq
    O\Parentheses{ 
        ( \epsilon^{2k} \vee \epsilon^k )
        ( 1 + \Norm{\mu}_\Ww + \Norm{\alpha^\epsilon}_\Aa + \Norm{\widehat{\alpha}}_\Aa )^{2k} 
    },
    \\
    \label{Eq:A-Priori-III}
    &\sup_{[0,T]} \Esp{ \Abs{ \phi_t }^{2k} }
    \leq
    O\Parentheses{ 
        \epsilon^k 
        ( 1 + \Norm{\mu}_\Ww + \Norm{\alpha^\epsilon}_\Aa + \Norm{\widehat{\alpha}}_\Aa )^{2k} 
    },
    \\
    \label{Eq:A-Priori-IV}
    &\sup_{[0,T]} \Esp{ \Abs{ \vartheta_t }^{2k} }
    \leq
    O\Parentheses{ 
        ( \epsilon^{2k} \vee \epsilon^{4k} )
        ( 1 + \Norm{\mu}_\Ww + \Norm{\alpha^\epsilon}_\Aa + \Norm{\widehat{\alpha}}_\Aa )^{4k} 
    }.
\end{align}

\end{lemma}

\begin{proof}
Under Assumptions \ref{Assumptions-C}, \eqref{Eq:A-Priori-I} is a direct result from the growth condition and \cite[Thm. 7.1]{mao_stochastic_2011}, as well as the elementary fact that
\begin{align*}
    &\Abs{ \sum_{i=1}^n x_i }^p \leq n^{p-1} \sum_{i=1}^n \Abs{ x_i }^p,
    &
    &p \geq 1,
\end{align*}
and Gronwall's inequality.

On the other hand, from the Lipschitz condition,
\begin{align*}
    &\sup_{[0,T]} \Esp{ \Abs{ \delta \widehat{b}_s^\mu( X^\epsilon, \alpha^\epsilon ) }^{k} }
    \leq 
    O\Parentheses{ \Parentheses{
        1 + \Norm{ \mu }_\Ww + \Norm{ \alpha^\epsilon }_\Aa + \Norm{ \widehat{\alpha} }_\Aa 
    }^{k} },
    \\
    &\sup_{[0,T]} \Esp{ 
        \Norm{ \delta \widehat{\sigma}_s^\mu( X^\epsilon, \alpha^\epsilon ) }^{2k}_\mathrm{Fr} 
    }
    \leq
    O\Parentheses{ \Parentheses{
        1 + \Norm{ \mu }_\Ww + \Norm{ \alpha^\epsilon }_\Aa + \Norm{ \widehat{\alpha} }_\Aa
    }^{2k} },
    \\
    &\sup_{[0,T]} \Esp{ 
        \Norm{ 
            \delta \widehat{\gamma}_s^\mu( X^\epsilon, \alpha^\epsilon) 
        }^{2k}_\mathrm{\lambda_t} 
    }
    \leq
    O\Parentheses{ \Parentheses{
        1 + \Norm{ \mu }_\Ww + \Norm{ \alpha^\epsilon }_\Aa + \Norm{ \widehat{\alpha} }_\Aa
    }^{2k} }.
\end{align*}
Thus, from the boundedness of the derivatives and \cite[Thm. 7.1]{mao_stochastic_2011},
\begin{align*}
    &\EE \Abs{ 
        \int_{ [0,t] \cap I^\epsilon } \delta \widehat{b}_s^\mu( X, \alpha ) \diff s 
    }^{2k}
    \leq
    O\Parentheses{ 
        \epsilon^{2k} 
        \Parentheses{ 
            1 + \Norm{ \mu }_\Ww + \Norm{ \alpha^\epsilon }_\Aa + \Norm{ \widehat{\alpha} }_\Aa
        }^{2k} 
    },
    \\
    &\EE \Abs{ 
        \int_{[0,t] \cap I^\epsilon } \delta \widehat{\sigma}_s^\mu( X, \alpha ) \diff W_s 
    }^{2k}
    \leq
    O\Parentheses{ 
        \epsilon^{k} 
        \Parentheses{
            1 + \Norm{ \mu }_\Ww + \Norm{ \alpha^\epsilon }_\Aa + \Norm{ \widehat{\alpha} }_\Aa
        }^{2k} 
    },
    \\
    &\EE \Abs{ 
        \int_{ ([0,t] \cap I^\epsilon) \times \mathbf{R} }
            \delta \widehat{\gamma}_s^\mu( X, \alpha \vert r )
        \tilde{\eta}^\lambda (\diff s, \diff r) 
    }^{2k}
    \leq
    O\Parentheses{ 
        \epsilon^{k} 
        \Parentheses{ 
            1 + \Norm{ \mu }_\Ww + \Norm{ \alpha^\epsilon }_\Aa + \Norm{ \widehat{\alpha} }_\Aa
        }^{2k} 
    },
\end{align*}
which in turn yields \eqref{Eq:A-Priori-II}-\eqref{Eq:A-Priori-IV}.

\end{proof}

\begin{proposition}
\label{Prop:Second-order-expanstion-of-spike-variation}
Given a random environment $\mu$, let $(X,\alpha)$ and $(\widehat{X},\widehat{\alpha})$ be as in Lemma \ref{Lemma:A-Priori-Estimates}. Then, for  $\epsilon \in (0,1)$,
\begin{align}
    \label{Eq:Second-order-expansion-X^epsilon}
    &\sup_{[0,T]} \Esp{ 
        \Abs{ \delta \widehat{X}^\epsilon_t - \phi_t - \vartheta_t }^{2} 
    }
    =
    o\Parentheses{ \epsilon^{2} }.
\end{align}
\end{proposition}

\begin{proof}

Using the integral remainder from the Taylor expansion \eqref{Eq:Taylor-expansion} on the coefficients, from \eqref{Eq:delta-hat-X^epsilon}-\eqref{Eq:vartheta} we obtain:
\begin{align*}
    \delta \widehat{X}^\epsilon_t - \phi_t - \vartheta_t
    =&
    X^\epsilon_t
    -
    \overbrace{ \Parentheses{ \widehat{X}_t + \phi_t + \vartheta_t } }^{ =: \rho_t }
    \\
    =&
    \int_0^t \Braces{
        \delta b_s^\mu\Parentheses{ X^\epsilon, \rho \middle\vert \alpha^\epsilon } + b^\star_s
    } \diff s
    +
    \int_0^t \Braces{
        \delta \sigma_s^\mu\Parentheses{ X^\epsilon, \rho \middle\vert \alpha^\epsilon } + \sigma^\star_s
    } \diff W_s
    \\
    &+
    \int_{ [0,t] \times \mathbf{R} } \Braces{
        \delta \gamma_s^\mu\Parentheses{ X^\epsilon, \rho \middle\vert \alpha^\epsilon, r } 
        + 
        \gamma^\star_{s-}(r)
    } \diff \tilde{\eta}^\lambda( \diff s, \diff r),
\end{align*}
where
\begin{align*}
    \Psi^\star_t(\cdot)
    :=&
    \delta \Diff \widehat{\Psi}_t^\mu( \alpha^\epsilon \vert \cdot, \widehat{X} ) \vartheta_t
    +
    \Parentheses{ \Id_d \otimes \vartheta_t \cdot \Diff^2 \widehat{ \Psi }_t^\mu(\cdot) \phi_t }
    +
    \Parentheses{ 
        \Id_d \otimes ( \phi_t + \vartheta_t ) \cdot \Diff^2 \widehat{ \Psi }_t^\mu(\cdot) \vartheta_t
    },
    \\
    &+
    \int_0^1 \Parentheses{
        \Id_d \otimes  (\phi_t + \vartheta_t)
        \cdot 
        \delta \Diff^2 \widehat{ \Psi }_t^\mu \Parentheses{ 
            \widehat{X} + (1 - \tau)(\phi + \vartheta), \alpha^\epsilon \middle\vert \cdot
        } (\phi_t + \vartheta_t)
    } \diff \tau,
\end{align*}
for $\Psi = b,\sigma,\gamma$. Moreover, from the integral remainder of Taylor's expansion, the identity
\begin{align}
    \label{Eq:Remainder-first-order-expansion}
    &\delta \Psi_t^\mu
    \Parentheses{ X^\epsilon, \rho \middle\vert \alpha^\epsilon }
    =
    \int_0^1 \Diff \Psi_t^\mu \Parentheses{
        \rho + (1 - \tau)(X^\epsilon - \rho), \alpha^\epsilon
    } (X^\epsilon_t - \rho_t) \diff \tau
\end{align}
holds for $\diff \PP \otimes \diff t$-almost all $(\omega, t) \in \Omega \times [0,T]$. Hence, from the boundedness of the derivatives,
\begin{align*}
    \Abs{ \delta \Psi_s^\mu \Parentheses{ X^\epsilon, \rho \middle\vert \alpha^\epsilon,\cdot} }^2
    &\leq
    C \Abs{ X^\epsilon_s - \rho_s }^2,
    &
    &\diff \PP \otimes \diff s-\mbox{a.e.}
\end{align*}

Applying Itô's formula and Gronwall's inequality yields
\begin{align}
    \label{Eq:Upper-bound-for-variation-of-X^epsilon}
    \sup_{[0,T]} \Esp{ \Abs{ X^\epsilon_t - \rho_t }^2 }
    \leq
    C \sup_{[0,T]} \EE \Bigg[  
        \Abs{ \int_0^t b^\star_s \diff s }^2
        &+
        \Abs{ \int_0^t \sigma^\star_s \diff W_s }^2
        \\ \nonumber
        &+
        \Abs{ \int_{[0,t] \times \mathbf{R}} 
            \gamma_{s-}^\star(r) \diff \tilde{\eta}^\lambda( \diff s, \diff r) }^2
    \Bigg].
\end{align}
However, from the estimates of Lemma \ref{Lemma:A-Priori-Estimates} and the boundedness of the second order derivatives, 
\begin{align*}
    \mbox{r.h.s. of \eqref{Eq:Upper-bound-for-variation-of-X^epsilon}}
    \leq
    O\Parentheses{
        \epsilon^3 
        ( 1 + \Norm{ \mu }_\Ww + \Norm{ \alpha^\epsilon }_\Aa + \Norm{ \widehat{\alpha} }_\Aa )^8 
    },
\end{align*}
thus completing the proof.

\end{proof}

Regarding the cost functional,  using the fact that $(\widehat{X},\widehat{\alpha})$ is optimal, it is clear that
\begin{align}
    \label{Eq:Inequality-of-optimal-pair}
    \delta \widehat{J}^\mu(X,\alpha)
    &=
    \Esp{ \int_0^T \delta \widehat{f}_t^\mu(X,\alpha) \diff t + \delta \widehat{g}_T^\mu(X) } 
    \geq 
    0
\end{align}
for all $(X,\alpha) \in \Ss^2(\RR^d) \times \Aa$. At the same time, from Assumptions \ref{Assumptions-C} and previous estimates we have that
\begin{align*}
    \Esp{ \Abs{
        \delta \widehat{g}_T^\mu(X^\epsilon) 
        -
        \Diff \widehat{ g }_T^\mu ( \delta \widehat{X}^\epsilon_T )
        -
        \frac{1}{2} \Parentheses{ 
            \delta \widehat{X}^\epsilon_T \cdot \Diff^2 \widehat{ g }_T^\mu ( \delta \widehat{X}^\epsilon_T )
        }
    }^2 } 
    \leq
    C \Esp{ \Abs{ \delta \widehat{X}^\epsilon_T }^6 }
    =
    O(\epsilon^3)
\end{align*}
for all $\epsilon > 0$ sufficiently small, with an analogue result for $f$. Our goal now is to define the adjoint variational system satisfied by the pay-off $J^\mu$ and derive the corresponding maximum principle.

\subsubsection{Necessary conditions through the Hamiltonian system}

\begin{theorem}
\label{Thm:Necessary-PMP}
Let $(\widehat{X},\widehat{\alpha})$ be an optimal pair for the optimization problem \eqref{Eq:MFG-Payoff}-\eqref{Eq:MFG-System} with random environment $\mu$. Then, there exists a pair of random processes 
\begin{align}
    \tag{\ref{Eq:YZUM}}
    (Y,Z,U,M) 
    \in
    \Ss^2(\RR^d) \times \HH^{2}(W) \times \HH^{2}(\tilde{\eta}^\lambda) \times \Hh^{2,\perp}
\end{align}
and
\begin{align}
    \tag{\ref{Eq:Tildes-YZUM}}
    \Parentheses{\widetilde{Y},\widetilde{Z},\widetilde{U},\widetilde{M}}
    \in
    \Ss^2(\RR^{d \times d}) \times \HH^{2}(W)^k \times \HH^{2}(\tilde{\eta}^\lambda)^l \times \Hh^{2,\perp}
\end{align}
such that 
\begin{enumerate}
    \item 
        $(Y,Z,U,M)$ is a solution to the first order adjoint BSDE \eqref{Eq:Y};

    \item
        $\Parentheses{\widetilde{Y},\widetilde{Z},\widetilde{U},\widetilde{M}}$ is a solution to the second-order adjoint BSDE \eqref{Eq:tilde-Y};

    \item
        \textbf{Extended Pontryagin condition:} for any admissible control $\alpha \in \Aa$, the inequality
        \begin{align}
            \label{Eq:Extended-Pontryagin-condition}
            \delta \widehat{H}_t^\mu \Parentheses{ \alpha \middle\vert \widehat{X},Y,Z,U }
            +&
            \sum_{i=1}^k \Parentheses{ 
                \delta \widehat{ \sigma }^{ \mu (\cdot,i) }_t (\alpha \vert \widehat{X})
                \cdot
                \frac{\widetilde{Y}_{s-}}{2} 
                    \delta \widehat{ \sigma }^{ \mu (\cdot,i) }_t (\alpha \vert \widehat{X})
            }_{\mathrm{Fr}}
            \\ \nonumber
            &+
            \sum_{j=1}^l \Parentheses{ 
                \delta\widehat{ \gamma }^{\mu (\cdot,j)}_t (\alpha \vert \widehat{X})
                \cdot
                \frac{\widetilde{Y}_{s-}}{2} 
                    \delta\widehat{ \gamma }^{\mu (\cdot,j)}_t (\alpha \vert \widehat{X}) 
            }_{\lambda_t^{(j)}}
            \geq
            0
        \end{align}
        holds $\diff \PP \otimes \diff t$-a.e. 
\end{enumerate}

\end{theorem}

\begin{remark}
In accordance with the   notation introduced in Section \ref{Sec:Basic-Notation}, the term $\delta \widehat{H}_t^\mu \big( \alpha \vert \widehat{X}, Y, Z$, $U \big)$ in \eqref{Eq:Extended-Pontryagin-condition} denotes the difference of Hamiltonians due to changes away from the optimal control, whenever the adjoint FBSDE is ``freezed'' at $(\widehat{X}_t$, $Y_t$, $Z_t$, $U_t$, $M_t)$.
\end{remark}

\begin{proof}
Observe that by having the arbitrary random environment $\mu$ and the optimal control $\widehat{\alpha}$ fixed, equations \eqref{Eq:MFG-System} and \eqref{Eq:Y}-\eqref{Eq:tilde-Y} form a FBSDE in environment-dependent jumps within an admissible set-up. Then, from Theorem 7 in \cite{hernandez-hernandez_coupled_2023}, there exists a Hamiltonian system $( \widehat{X},\widehat{Y},\widehat{Z},\widehat{U},\widehat{M} )$ as in \eqref{Eq:Hats-XYZUM}. Hence, we only need to prove (3).

Let $0 < \epsilon \ll 1$.  From the estimates computed at Lemma \ref{Lemma:A-Priori-Estimates} and Proposition \ref{Prop:Second-order-expanstion-of-spike-variation}, 
\begin{align*}
    &\Esp{ \delta \widehat{f}_t^\mu(X^\epsilon,\alpha^\epsilon) }
    =
    \Esp{ \delta \widehat{f}_t^\mu(\rho,\alpha^\epsilon) } + o(\epsilon^2)
    &
    &\mbox{and}
    &
    &\Esp{ \delta \widehat{g}_T^\mu(X^\epsilon) }
    =
    \Esp{ \delta \widehat{g}_T^\mu(\rho) } + o(\epsilon^2),
\end{align*}
where $\rho = \widehat{X} + \phi + \vartheta$. Then,
\begin{align}
    \label{Eq:Expansion-of-payoff-inequality}
    \delta\widehat{J}^\mu &(\alpha^\epsilon,X^\epsilon)
    =
    \Esp{ 
        \int_0^T \delta \widehat{f}_t^\mu( \rho, \alpha^\epsilon ) \diff t
        +
        \delta \widehat{g}_T^\mu(\rho)
    }
    +
    o(\epsilon^2)
    \\ \nonumber
    =&
    \Esp{ \int_0^T 
        \delta f_t^\mu \Parentheses{ \alpha^\epsilon, \widehat{\alpha} \middle\vert \rho } 
    \diff t }
    +
    \Esp{ 
        \int_0^T \delta \widehat{f}_t^\mu( \rho \vert \widehat{\alpha} ) \diff t
        +
        \delta \widehat{g}_T^\mu(\rho)
    }
    +
    o(\epsilon^2)
    \\ \nonumber
    =&
    \Esp{ \int_0^T 
        \delta f_t^\mu \Parentheses{ \alpha^\epsilon, \widehat{\alpha} \middle\vert \rho } 
    \diff t }
    +
    \Esp{
        \int_0^T \Diff \widehat{ f }_t^\mu \Parentheses{ \phi_t + \vartheta_t } \diff t
        +
        \Diff \widehat{ g }_T^\mu \Parentheses{ \phi_T + \vartheta_T } 
    }
    \\ \nonumber
    &+
    \frac{1}{2} \Esp{
        \int_0^T \Parentheses{ 
            ( \phi_t + \vartheta_t ) \cdot \Diff^2 \widehat{ f }_t^\mu ( \phi_t + \vartheta_t ) 
        } \diff t
        +
        \Parentheses{ 
            ( \phi_T + \vartheta_T ) \cdot \Diff^2 \widehat{ g }_T^\mu ( \phi_T + \vartheta_T ) 
        }
    }
    \\ \nonumber
    &+
    o(\epsilon^2).
\end{align}

Now, from equations \eqref{Eq:Diff_H}-\eqref{Eq:Diff^2_H} and \cite[Lemma 15]{hernandez-hernandez_coupled_2023} we obtain the following identity for the first order variation:
\begin{align}
    \label{Eq:Payoff-first-order-variation}
    &\Esp{ 
        \int_0^T \Diff \widehat{ f }_t^\mu \Parentheses{ \phi_t + \vartheta_t } \diff t
        +
        \Diff \widehat{ g }_T^\mu \Parentheses{ \phi_T + \vartheta_T } 
    }
    \\ \nonumber
    &=
    \Esp{ \Parentheses{ Y_0 \cdot \delta \widehat{X}_0^\epsilon  } }
    +
    \Esp{ \int_0^T \Braces {
        \delta \widehat{H}_t^\mu(\alpha^\epsilon \vert \widehat{X},Y,Z,U)
        -
        \delta \widehat{f}_t^\mu(\alpha^\epsilon \vert \widehat{X} )
    } \diff t }
    \\ \nonumber
    &\quad-
    \sum_{i=1}^k \Esp{ \int_0^T
        \Parentheses{ 
            Z_t^{(i,\cdot)} 
            \cdot 
            \delta \Diff \widehat{ \sigma }_t^{\mu(\cdot,i)}( \alpha^\epsilon \vert \widehat{X} ) \phi_t 
        }
    \diff t }
    \\ \nonumber
    &\quad-
    \sum_{j=1}^l \Esp{ \int_0^T
        \Parentheses{ 
            U_{t-}^{(j,\cdot)}
            \cdot
            \delta \Diff \widehat{ \gamma }_{t-}^{\mu(\cdot,j)}( \alpha^\epsilon \vert \widehat{X} ) \phi_{t-}
        }_{ \lambda_t^{(j)} }
    \diff t }
    \\ \nonumber
    &\quad-
    \frac{1}{2} \Esp{ \int_0^T 
        \Parentheses{
            \phi_t
            \cdot 
            \Braces{ \Diff^2 \widehat{ H }_t^\mu(Y,Z,U) - \Diff^2 \widehat{ f }_t^\mu } \phi_t 
        }
    \diff t }.
\end{align}

However, because of the assumptions of boundedness for the derivatives and the square integrability of $(Y,Z,U,M)$, the \emph{a priori} estimates computed in Lemma \ref{Lemma:A-Priori-Estimates} yield
\begin{align*}
    &\Abs{ \Esp{ 
        \Parentheses{ Y_0 \cdot \delta \widehat{X}_0^\epsilon } 
        +
        \frac{1}{2} \int_0^T \Parentheses{
            \phi_t
            \cdot 
            \Braces{ \Diff^2 \widehat{ H }_t^\mu(Y,Z,U) - \Diff^2 \widehat{ f }_t^\mu } \phi_t 
        } \diff t 
    } }
    \leq
    o(\epsilon).
\end{align*}
Similarly, for all $1 \leq i \leq k$ and all $1 \leq j \leq l$:
\begin{align*}
    &\Abs{ \Esp{ \int_0^T
        \Parentheses{ 
            Z_t^{(i,\cdot)} 
            \cdot 
            \delta \Diff \widehat{ \sigma }_t^{\mu (\cdot,i)}( \alpha^\epsilon \vert \widehat{X} ) \phi_t 
        }
    \diff t } }
    \\
    &{\hskip0.1\linewidth}+
    \Abs{ \Esp{ \int_0^T
        \Parentheses{ 
            U_t^{(j,\cdot)}
            \cdot
            \delta \Diff \widehat{ \gamma }_t^{\mu(\cdot,j)}( \alpha^\epsilon \vert \widehat{X} ) \phi_{t-}
        }_{ \lambda_t^{(j)} }
    \diff t } }
    \leq
    o(\epsilon).
\end{align*}
Then, plugging \eqref{Eq:Payoff-first-order-variation} into \eqref{Eq:Expansion-of-payoff-inequality}:
\begin{align}
    \label{Eq:Expansion-of-payoff-inequality'}
    \tag{\ref{Eq:Expansion-of-payoff-inequality}'}
    \delta\widehat{J}^\mu &(X^\epsilon,\alpha^\epsilon)
    =
    \Esp{ \int_0^T \Braces {
        \delta f_t^\mu \Parentheses{ \alpha^\epsilon, \widehat{\alpha} \middle\vert \rho } 
        -
        \delta \widehat{ f }_t^\mu(\alpha^\epsilon \vert \widehat{X} )
    } \diff t }
    \\ \nonumber
    &+
    \frac{1}{2} \Esp{ 
        \int_0^T \Parentheses{ 
            ( \phi_t + \vartheta_t ) \cdot \Diff^2 \widehat{ f }_t^\mu ( \phi_t + \vartheta_t ) 
        } \diff t
        +
        \Parentheses{ 
            ( \phi_T + \vartheta_T ) \cdot \Diff^2 \widehat{ g }_T^\mu ( \phi_T + \vartheta_T ) 
        } 
    }
    \\ \nonumber
    &+
    \Esp{ \int_0^T \delta \widehat{ H }_t^\mu(\alpha^\epsilon \vert \widehat{X},Y,Z,U) \diff t }
    +
    o(\epsilon).
\end{align}

Observe that for the first term at the r.h.s. of \eqref{Eq:Expansion-of-payoff-inequality'},
\begin{align}
    \label{Eq:Variation-on-tilde-f}
    \EE & \Brackets{ \int_0^T \Braces {
        \delta f_t^\mu \Parentheses{ \alpha^\epsilon, \widehat{\alpha} \middle\vert \rho } 
        -
        \delta \widehat{f}_t^\mu(\alpha^\epsilon \vert \widehat{X} )
    } \diff t }
    \\ \nonumber
    =&
    \Esp{ \int_0^T \Braces{
        \delta f_t^\mu \Parentheses{ \rho, \widehat{X} \middle\vert \alpha^\epsilon } 
        -
        \delta \widehat{f}_t^\mu( \rho \vert \widehat{\alpha} )
    } \diff t } 
    \\ \nonumber
    =&
    \EE \bigg[  \int_0^T \Big\lbrace
        \delta \Diff \widehat{ f }_t^\mu \Parentheses{ \alpha^\epsilon \middle\vert \widehat{X} } 
            ( \phi_t + \vartheta_t ) 
        \\ \nonumber
        &{\hskip0.1\linewidth}+
        \frac{1}{2} \Parentheses{
            ( \phi_t + \vartheta_t ) 
            \cdot 
            \delta \Diff^2 \widehat{ f }_t^\mu \Parentheses{ \alpha^\epsilon \middle\vert \widehat{X} } 
                ( \phi_t + \vartheta_t ) 
        } 
    \Big\rbrace \diff t \bigg]
    +
    o(\epsilon^2).
\end{align}

Likewise, by replicating the proof of Lemma 15 in \cite{hernandez-hernandez_coupled_2023} twice, first on the product $\Diff^2 \widehat{ g }_T^\mu ( \phi_T + \vartheta_T )$ and then on the quadratic form $\Parentheses{ ( \phi_T + \vartheta_T ) \cdot \Diff^2 \widehat{ g }_T^\mu ( \phi_T + \vartheta_T ) }$, we obtain that the following identity for the second order variation:
\begin{align*}
    &\EE \Brackets{ \Parentheses{ 
        ( \phi_T + \vartheta_T )
        \cdot 
        \Diff^2 \widehat{ g }_T^\mu ( \phi_T + \vartheta_T ) 
    } }
    =
    \EE \Brackets{ \Parentheses{ 
        \delta \widehat{X}^\epsilon_0 \cdot \widetilde{Y}_0\ \delta \widehat{X}^\epsilon_0
    } }
    \\ \nonumber
    &+
    \int_0^T \EE \bigg[
        \Big( 
            (\phi_t + \vartheta_t ) \cdot (\tilde{Y}_t - B^{\star\star}_t) (\phi_t + \vartheta_t ) 
        \Big)
        +
        \Big( 
            (\phi_t + \vartheta_t ) \cdot \Brackets{ \tilde{Y}_t + \tilde{Y}_t^\top } b_t^{\star\star} 
        \Big)
    \bigg] \diff t 
    \\ \nonumber
    &-
    \sum_{i=1}^j \int_0^T \EE \bigg[
        \Parentheses{
            ( \phi_t + \vartheta_t ) \cdot \Sigma^{\star\star(i)}_t ( \phi_t + \vartheta_t )
        }
        -
        \Parentheses{
            (\phi_t + \vartheta_t )
            \cdot 
            \Brackets{ \tilde{Z}^{(i)}_t + \tilde{Z}^{(i)\top}_t } \sigma_t^{\star\star(\cdot,i)}
        }
        \\ \nonumber
        &{\hskip0.16\linewidth}-
        \Parentheses{ \sigma_t^{\star\star(\cdot,i)} \cdot \tilde{Y}_t \sigma_t^{\star\star(\cdot,i)} }
    \bigg]  \diff t
    \\ \nonumber
    &-
    \sum_{j=1}^l \int_0^T \EE \bigg[
        \Parentheses{
            ( \phi_t + \vartheta_t ) \cdot \Gamma^{\star\star(j)}_t ( \phi_t + \vartheta_t )
        }_{ \lambda^{(j)}_t }
        -
        \Parentheses{
            (\phi_t + \vartheta_t )
            \cdot 
            \Brackets{ \tilde{U}^{(j)}_t + \tilde{U}^{(j)\top}_t } \gamma_t^{\star\star(\cdot,j)}
        }_{ \lambda^{(j)}_t }
        \\ \nonumber
        &{\hskip0.16\linewidth}-
        \Parentheses{ 
            \gamma_t^{\star\star(\cdot,j)}
            \cdot 
            \tilde{Y}_t\gamma_t^{\star\star(\cdot,j)}
        }_{ \lambda^{(j)}_t }
    \bigg] \diff t,
\end{align*}
where $B^{\star\star}$, $\Sigma^{\star\star}$ and $\Gamma^{\star\star}$ are defined as in \eqref{Eq:B^starstar}, \eqref{Eq:Sigma^starstar} and \eqref{Eq:Gamma^starstar}, and $b^{\star\star}$, $\sigma^{\star\star}$ and $\gamma^{\star\star}$ denote the drift, diffusion and jump coefficients of the process $\phi + \theta$, respectively.

Because of the boundedness conditions on the derivatives and the square-integrability of the process $(\widetilde{Y},\widetilde{Z},\widetilde{U},\widetilde{M})$,  the linear and quadratic terms on $\phi + \vartheta$ are of order $o(\epsilon^2)$:
\begin{align}
    \label{Eq:Payoff-second-order-variation}
    &\Esp{ \Parentheses{ 
        (\phi_T + \vartheta_T) \cdot \Diff^2 \widehat{ g }_T^\mu (\phi_T + \vartheta_T) 
    } }
    \\ \nonumber
    &=
    \int_0^T \Esp{ 
        \Parentheses{ 
            \delta \widehat{ \sigma }_t^\mu (\alpha^\epsilon \vert \widehat{X} ) 
            \cdot 
            \widetilde{Y}_t\ \delta \widehat{ \sigma }_t^\mu (\alpha^\epsilon \vert \widehat{X} )
        }_{\mathrm{Fr}}
        +
        \Parentheses{ 
            \delta \widehat{ \gamma }_t^\mu (\alpha^\epsilon \vert \widehat{X} ) 
            \cdot 
            \widetilde{Y}_t\ \delta \widehat{ \gamma }_t^\mu (\alpha^\epsilon \vert \widehat{X} )
        }_{ \lambda_t }
    } \diff t
    +
    o(\epsilon^2).
\end{align}

Then, from \eqref{Eq:Expansion-of-payoff-inequality'}, \eqref{Eq:Variation-on-tilde-f} and \eqref{Eq:Payoff-second-order-variation}:
\begin{align}
    \nonumber
    \frac{ \delta\widehat{J}^\mu(X^\epsilon,\alpha^\epsilon) }{\Abs{\epsilon}}
    =&
    \frac{1}{\Abs{\epsilon}} \int_0^T \EE \bigg[  
        \delta \widehat{H}_t^\mu(\alpha^\epsilon\vert\widehat{X},Y,Z,U)
        +
        \Big(
            \delta \widehat{ \sigma }_t^\mu (\alpha^\epsilon \vert \widehat{X} ) 
            \cdot 
            \frac{\widetilde{Y}_t}{2}\ \delta \widehat{ \sigma }_t^\mu (\alpha^\epsilon \vert \widehat{X} )
        \Big)_{\mathrm{Fr}}
    \bigg] \diff t
    \\ \label{Eq:Almost-derivative-of-Payoff}
    &+
    \frac{1}{ \Abs{\epsilon} } \int_0^T \EE \bigg[  
        \Big(
            \delta \widehat{ \gamma }_t^\mu (\alpha^\epsilon \vert \widehat{X} ) 
            \cdot 
            \frac{\widetilde{Y}_t}{2}\ \delta \widehat{ \gamma }_t^\mu (\alpha^\epsilon \vert \widehat{X} )
        \Big)_{ \lambda_t }
    \bigg] \diff t
    +
    \frac{o(\Abs{\epsilon})}{ \Abs{\epsilon} }.
\end{align}

To guarantee the differentiation with respect to $\epsilon$, define $l:[0,T] \times A \to \RR$ as
\begin{align*}
    l(t,a)
    :=
    \EE \bigg[  
        &\delta \widehat{H}_t^\mu( a \vert\widehat{X},Y,Z,U)
        +
        \Big(
            \delta \widehat{ \sigma }_t^\mu ( a \vert \widehat{X} ) 
            \cdot 
            \frac{\widetilde{Y}_t}{2}\ \delta \widehat{ \sigma }_t^\mu ( a \vert \widehat{X} )
        \Big)_{\mathrm{Fr}}
        \\
        &+
        \Big(
            \delta \widehat{ \gamma }_t^\mu ( a \vert \widehat{X} ) 
            \cdot 
            \frac{\widetilde{Y}_t}{2}\ \delta \widehat{ \gamma }_t^\mu ( a \vert \widehat{X} )
        \Big)_{ \lambda_t }
    \bigg].
\end{align*}
From \eqref{Eq:alpha^epsilon}, \eqref{Eq:Inequality-of-optimal-pair} and \eqref{Eq:Almost-derivative-of-Payoff}, an application of Lyapunov's convexity theorem \cite[Lemma 1]{kappel_maximum_1985} yields
\begin{align}
    \label{Eq:Derivative-of-Payoff}
    \frac{ \diff \widehat{J}^\mu(X^\alpha,\alpha) }{ \diff \alpha } \Bigg\vert_{\alpha = \widehat{\alpha}}
    =
    \lim_{ \Abs{\epsilon} \to 0^+ } \frac{1}{\Abs{\epsilon}} \int_{I_\epsilon} l(t, \alpha) \diff t
    =
    \int_0^T l(t, \alpha) \diff t
    \geq 0
\end{align}
$\PP$-a.s. for any arbitrary admissible control $\alpha \in \Aa$. Using the same procedure as in \cite{kushner_necessary_1972} (see also \cite{xun_yu_unified_1991} or \cite{tang_necessary_1994}), \eqref{Eq:Derivative-of-Payoff} implies \eqref{Eq:Extended-Pontryagin-condition}, thus concluding the proof.
\end{proof}

\subsection{Sufficient conditions for optimality}

We now complement the results  given in  Theorem \ref{Thm:Necessary-PMP} by presenting a converse result describing a set of  sufficient conditions that a minimizer must satisfy. Our approach is based on \cite{framstad_sufficient_2004}.

\begin{theorem}
\label{Thm:Sufficent-PMP}
Let $(\widehat{X},\widehat{\alpha}) \in \Ss^2(\RR^d) \times \Aa$ be an $\widehat{\alpha}$-controlled pair.   If
\begin{itemize}

    \item
        $g$ is a convex function on $x$;
        
    \item 
        there exist a pair of processes $(Y,Z,U,M)$ and $(\widetilde{Y},\widetilde{Z},\widetilde{U},\widetilde{M})$, defined as in \eqref{Eq:YZUM} and \eqref{Eq:Tildes-YZUM}, that solve the first and second order adjoint BSDEs \eqref{Eq:Y} and \eqref{Eq:tilde-Y}, respectively;

    \item
        the mapping
        \begin{align}
            \label{Eq:Hh}
            x \longmapsto \min_{a \in A} H_t^\mu(x,Y,Z,U,a) =: \Hh_t^\mu(x),
        \end{align}
         which represents the minimized Hamiltonian at the state $x$, exists and is convex $\diff \PP \otimes \diff t$-a.e.;

    \item
        for any admissible control $\alpha \in \Aa$, the inequality
        \begin{align}
            \label{Eq:Pontryagin-condition}
            &\delta \widehat{H}_t^\mu \Parentheses{ \alpha \middle\vert \widehat{X}, Y,Z,U } \geq 0
        \end{align}
        holds $\diff \PP \otimes \diff t$-a.a.;
        
\end{itemize}
then, $(\widehat{X},\widehat{\alpha})$ is an optimal pair for the control problem in random environment $\mu$ given by \eqref{Eq:MFG-Payoff}-\eqref{Eq:MFG-System}.
\end{theorem}

\begin{proof} 
Let $\alpha \in \Aa$ be an arbitrary control, and set $X = X^\alpha$. In order to show that \eqref{Eq:Inequality-of-optimal-pair} holds, recall that
\begin{align*}
    \delta\widehat{J}^\mu( X, \alpha )
    =
    \Esp{ \int_0^T \delta \widehat{f}_t^\mu(X, \alpha) \diff t + \delta \widehat{g}_T^\mu(X) }.
\end{align*}
From the convexity of $g$ and \cite[Lemma 15]{hernandez-hernandez_coupled_2023},
\begin{align*}
    \Esp{ \delta \widehat{g}_T^{\mu}(X) }
    \geq& 
    \Esp{ \Parentheses{ \delta \widehat{X}_T \cdot \Diff \widehat{ g }_T^{\mu\top} } }
    \\
    =& 
    \int_0^T \Esp{
        \Parentheses{ \delta \widehat{b}_t^\mu(X,\alpha) \cdot Y_t }
        +
		\Parentheses{ \delta \widehat{\sigma}_t^\mu(X,\alpha) \cdot Z_t }_{ \mathrm{Fr} }
    } \diff t
    \\
    &+
    \int_0^T \Esp{
        \Parentheses{ \delta \widehat{\gamma}_t^\mu(X,\alpha) \cdot U_{t} }_{ \lambda_t }
        -
        \Parentheses{ \delta \widehat{X}_t \cdot \Diff \widehat{ H }_t^{\mu\top}(Y,Z,U) }
    } \diff s.
\end{align*}
Hence, from \eqref{Eq:H}, \eqref{Eq:Hh} and the convexity of $\Hh^\mu$,
\begin{align*}
    \delta\widehat{J}^\mu( X, \alpha )
    \geq&
    \int_0^T \Esp{
        \delta \widehat{H}_t^\mu( X, \alpha \vert Y,Z,U)
        -
        \Parentheses{ \delta \widehat{X}_t \cdot \Diff \widehat{ H }_t^{\mu\top}(Y,Z,U) }
    } \diff t
    \\
    \geq&
    \int_0^T \Esp{
        \delta \widehat{H}_t^\mu( X, \alpha \vert Y,Z,U)
        -
        \delta \widehat{\Hh}_t^\mu( X )
    } \diff t.
\end{align*}

To conclude, observe that  \eqref{Eq:Pontryagin-condition} implies that
\begin{align*}
    H_t^\mu( \widehat{X}, Y, Z, U, \widehat{\alpha} )
    =
    \min_{a \in A} H_t^\mu( \widehat{X}, Y, Z, U, a ),
\end{align*}
and from \eqref{Eq:Hh},
\begin{align*}
    H_t^\mu( x, Y, Z, U, \alpha )
    &\geq
    \Hh_t^\mu(x),
    &
    &\forall x \in \RR^d,
\end{align*}
$\PP$-a.s for all $t \in [0,T]$. As a result,
\begin{align*}
    \delta \widehat{H}_t^\mu ( x, \alpha \vert Y, Z, U )
    &\geq
    \delta \widehat{\Hh}_t^\mu ( x ),
    &
    &\forall x \in \RR^d,
\end{align*}
therefore completing the proof.
\end{proof}

\section{Solution to the Mean-Field Game}
\label{Sec:Sol MFG}

Having described the elements to characterize an optimal solution when a fixed random environment $\mu$ is fixed, through Theorems \ref{Thm:Necessary-PMP} and \ref{Thm:Sufficent-PMP}, we go back to the Mean-Field Game described in Definition \ref{Def:MFE} in order to analyze the existence of an equilibrium.

Let $X$ be the solution to a conditional McKean-Vlasov equation given the information of the common noise process $(N^\lambda, \xi)$. Since the solution to a (conditional) McKean-Vlasov equation induces an \textit{admissible} set-up  \cite[Lemma 11]{hernandez-hernandez_conditional_2023}, we can proceed to apply the previous results to the \textit{admissible} random environment $\Ll^1(X)$. This is the key idea behind the next theorem.

\begin{theorem}
\label{Thm:Characterisation-of-Mean-Field-equilibrium}

\hspace{0.99\linewidth}\newline
\noindent
\textbf{(a) Necessary conditions for equilibrium.} Let $(\widehat{\mu},\widehat{X})$ be a strong MFE with a corresponding admissible control $\widehat{\alpha} \in \Aa$, such that 
\begin{align}
    \label{Eq:8th-moment-condition}
    \widehat{X}_t \in L^8(\PP)\mbox{ uniformly on }t \in [0,T].
\end{align}

Then, for $\diff \PP \otimes \diff t$-a.a $(\omega,t) \in \Omega \times [0,T]$ there exists a process $\big( \widehat{X}$, $\widehat{Y}$, $\widehat{Z}$, $\widehat{U}$, $\widehat{M} \big)$, defined as in \eqref{Eq:Hats-XYZUM} with $\widehat{M} \equiv 0$, that solves the Conditional McKean-Vlasov FBSDE
\begin{align}
    \label{Eq:McK-V-X}
    \widehat{X}_t
    =&
    \widehat{X}_0
    +
    \int_0^t \widehat{b}_s \diff s
    +
    \int_0^t \widehat{\sigma}_s \diff W_s
    +
    \int_{[0,t] \times \mathbf{R} } \widehat{\gamma}_{s}(r) \tilde{\eta}^\lambda (\diff s, \diff r),
    \\[8pt] \label{Eq:McK-V-Y}
    Y_t
    =&
    \Diff \widehat{ g }^\top_T
    +
    \int_t^T \Diff \widehat{ H }_s^\top \diff s
    -
    \int_t^T Z_s \diff W_s
    -
    \int_{(t,T] \times \mathbf{R} } U_{s}(r) \tilde{\eta}^\lambda (\diff s, \diff r),
    \\[8pt] \label{Eq:McK-V-Tilde-Y}
    \widetilde{Y}_t
    =&
    \Diff^2 \widehat{ g }_T^\mu
    +
    \int_t^T \widehat{B}^{\star\star}_s \diff s
    -
    \sum_{i=1}^k \int_t^T \Big\lbrace 
        \widetilde{Z}^{(i)}_s \diff W^{(i)}_s - \widehat{\Sigma}_s^{\star\star(i)} \diff s 
    \Big\rbrace
    \\ \nonumber
    &-
    \sum_{j=1}^l \int_{ (t,T] \times \mathbf{R} } \Big\lbrace
        \widetilde{U}^{(j)}_s(r) \tilde{\eta}^{\lambda(j)}( \diff s, \diff r )
        -
        \widehat{\Gamma}_s^{\star\star(j)}(r) \big(\diff s \otimes K^{\lambda(j)}_s (\diff r) \big)
    \Big\rbrace,
\end{align}
where $\widehat{B}^{\star\star}$, $\widehat{\Sigma}^{\star\star}$ and $\widehat{\Gamma}^{\star\star}$ are defined as in \eqref{Eq:B^starstar}, \eqref{Eq:Sigma^starstar} and \eqref{Eq:Gamma^starstar}, with $\mu = \widehat{\mu}$. In other words, $\big( \widehat{X}$, $\widehat{Y}$, $\widehat{Z}$, $\widehat{U}$, $\widehat{M} \big)$ is the Hamiltonian system for the optimization problem \eqref{Eq:MFG-Payoff}-\eqref{Eq:MFG-System}, under the endogenous random environment $\widehat{\mu} = \Ll^1(\widehat{X})$.

\noindent
\textbf{(b) Sufficient conditions for equilibrium.} Suppose that  the following conditions are satisfied:
\begin{itemize}

    \item
        $g$ is a convex function on $x$;
        
    \item 
        the mapping 
        \begin{align*}
            (\nu,x,y,z,u) 
            \longmapsto 
            \widehat{a} \in \argmin_{a \in A} H_t^\nu \Parentheses{ x, y, z, u, a }
        \end{align*}
        exists $\diff \PP \otimes \diff t$-almost everywhere;

    \item
        there exists a process $(\widehat{X}, \widehat{Y}, \widehat{Z}, \widehat{U}, \widehat{M})$ that solves the conditional McKean-Vlasov FBSDE \eqref{Eq:McK-V-X}-\eqref{Eq:McK-V-Tilde-Y} for the (admissible) control
        \begin{align*}
            \widehat{\alpha}_t
            :=
            \argmin_{a \in A} H_t^{\Ll^1(\widehat{X})} \Parentheses{ \widehat{X}, Y, Z, U, a };
        \end{align*}

    \item 
        $\widehat{X}$ satisfies the 8-th moment condition \eqref{Eq:8th-moment-condition};

    \item
        the mapping
        \begin{align*}
            x \longmapsto \min_{a \in A} H_t^{\Ll^1(\widehat{X})}( x, Y, Z, U, a) 
        \end{align*}
        exists and is convex $\diff \PP \otimes \diff t$-a.a.

\end{itemize}
Then, $\Parentheses{\Ll^1(\widehat{X}),\widehat{X}}$ is a strong MFE for the game \eqref{Eq:MFG-Payoff}-\eqref{Eq:MFG-System}. 

\end{theorem}

\begin{remark}
From our assumptions on $(\Omega,\FF)$ and the compatibility condition on admissible set-ups (condition 7 in Definition \ref{Def:Admissible-Set-Up}), there exists a regular conditional probability distribution
\begin{align*}
    \PP^1[\ \cdot\ ] := \Prob{\ \cdot\ \middle\vert \FF^{(X_0,\mu,\lambda,(N^\lambda,\xi))} },
\end{align*}
see \cite[Ch. 10.6]{bogachev_measure_2007} and \cite[Prop. 1.10]{carmona_probabilistic_2018-1} (alternatively, see the proof of Lemma 11 at \cite{hernandez-hernandez_conditional_2023}). In addition, we will also use $\EE^1$ (resp. $\operatorname{Var}^1$) to represent the expectation (resp. variance) operator with respect to $\PP^1$.
\end{remark}

\begin{proof}
From the properties of Wasserstein distances and the disintegration of measure $\PP$ with respect to $\PP^1$, 
\begin{align*}
    &\Esp{ W_{2}( \Ll^1(\widehat{X})_t, \delta_{ \{ 0 \} } )^8 } 
    \leq
    \Esp{ W_{8}( \Ll^1(\widehat{X})_t, \delta_{ \{ 0 \} } )^8 } 
    \leq 
    \Esp{ \Esp[^1]{ \Abs{ \widehat{X}_t }^{8} } }
\end{align*}
for all $t \in [0,T]$; then, from the moment condition at \eqref{Eq:8th-moment-condition} the space conditions in Assumptions \ref{Assumptions-C} are satisfied, while the rest are due to the coefficients and the controls themselves.

Now, for the first part, let $(\widehat{X}, \widehat{\alpha})$ be the optimal pair for the optimization in random environment \eqref{Eq:MFG-Payoff}-\eqref{Eq:MFG-System}, with $\mu = \Ll^1(\widehat{X})$ --which exists in view of the assumption that $( \widehat{\mu}, \widehat{X} )$ is an equilibrium--. The conclusion follows from a direct application of Theorem \ref{Thm:Necessary-PMP}. The second part follows directly from Theorem \ref{Thm:Sufficent-PMP} and Definition \ref{Def:MFE}.
\end{proof}

\begin{remark}
\label{Cor:Uncontrolled-jump-diffusion}
When the diffusion and jump terms ($\sigma$ and $\gamma$ respecively) do not deppend on the control, then the Hamiltonian system $\big( \widehat{X}, \widehat{Y}, \widehat{Z}, \widehat{U}, \widehat{M} \big)$ from Theorem \ref{Thm:Characterisation-of-Mean-Field-equilibrium} is reduced to $\big(  \widehat{X}, Y, Z, U, M \big)$, where $\widehat{X} \in \Ss^2(\RR^{d})$ satisfies \eqref{Eq:McK-V-X} and $\Parentheses{ Y, Z, U, M }$, defined as in \eqref{Eq:YZUM} with $M=0$, is a solution to \eqref{Eq:McK-V-Y}. That is, $(\widetilde{Y}, \widetilde{Z}, \widetilde{U}, \widetilde{M})$ is identically zero. Furthermore, the extended Pontryagin condition \eqref{Eq:Extended-Pontryagin-condition} is reduced to \eqref{Eq:Pontryagin-condition}:
\begin{align}
    \tag{\ref{Eq:Pontryagin-condition}'}
    \label{Eq:Pontryagin-condition'}
    &\widehat{H}_t
    =
    H^{\widehat{\mu}}_t
        \Parentheses{ \widehat{X}, \widehat{Y}, \widehat{Z}, \widehat{U}, \widehat{\alpha} }
    =
    \operatornamewithlimits{ess\ inf}_{\alpha \in \Aa} H^{\widehat{\mu}}_t
        \Parentheses{ \widehat{X}, \widehat{Y}, \widehat{Z}, \widehat{U}, \alpha },
    &
    &\forall t \in [0,T].
\end{align}
\end{remark}

Observe that the results we have presented so far are consistent with their analogous  analyzed in \cite{carmona_probabilistic_2018}-\cite{carmona_probabilistic_2018-1} and the references there-in; namely, those regarding the solution of MFGs with a common Gaussian noise.

\section{Examples and applications}
\label{Sec:EyA}

In this section we apply the  results obtained in the previous sections  to solve a \textit{Linear-Quadratic Mean-Field Game} arising from dynamical systems with self-exciting jumps, similar to the one presented in Example 19 from \cite{hernandez-hernandez_coupled_2023}. Furthermore,
in order  to illustrate the scope  of Theorems \ref{Thm:Necessary-PMP}-\ref{Thm:Sufficent-PMP} beyond MFGs, we address the problem of \textit{Mean-Variance Portafolio Selection within a Regime Switching Market} in Section \ref{Sec:RS-Porfolio-Selection}. This problem is inspired by Example 4.5 in \cite{nguyen_general_2021}, and makes use of the existence and convergence results shown previously in Example 19 and Section 5 from \cite{hernandez-hernandez_coupled_2023} and \cite{hernandez-hernandez_conditional_2023}, respectively.


\subsection{Linear-quadratic mean-field game}
\label{Sec:LQMFG}

Consider the  one dimensional MFG with payoff 
\begin{align}
    \nonumber
    J^\mu(X,\alpha)
    :=
    \EE\Bigg[
		&\int_0^T \Braces{
			\frac{f_1}{2} \Abs{ X_t - \Esp[^1]{X_t} }^2
			+
			f \Parentheses{ X_t - \Esp[^1]{X_t} } \alpha_t
			+
			\frac{f_2}{2} \Abs{ \alpha_t }^2
		} \diff t 
		\\ \label{Eq:(LQMFG)-J}
		&+
		\frac{g}{2} \Abs{ X_T - \Esp[^1]{X_T} }^2
	\Bigg],
\end{align}
where $f, f_1, f_2, g$ real constants satisfying 
\begin{align}
    \label{Eq:(LQFMG)Coefficient-conditions}
    &\Abs{b} f_2 = 2
    &
    &\mbox{and}
    &
    &f_1 f_2 > \frac{f^2}{2} + 1,
\end{align}

The state $X$ evolves according  to the controlled dynamics
\begin{align}
    \label{Eq:(LQMFG)-X}
    X_t 
    =&
    X_0
    +
    \int_0^t \Braces{ b \Parentheses{ X_{s} - \Esp[^1]{X_s} } + \alpha_s } \diff s
    +
    \sigma W_t
    +
    \gamma \tilde{N}^\lambda_t,
\end{align}
such that $b,\sigma,\gamma \in \RR$, and $\lambda$ is the additive intensity candidate \cite{hernandez-hernandez_coupled_2023}:
\begin{align}
    \label{Eq:LQ-lambda}
    \lambda_t 
    = 
    \lambda_0 + \psi^1 v_r(\mu_{t-}) + \int_{[0,t)} \psi^2(t-s) \diff N^\lambda_s,
\end{align}
where $\lambda_0$ and $\psi^1$ are positive constants, $\psi^2$ is a non-negative continuous deterministic function, and $v_r$ is the $\Pp_2(\RR)$-valued function defined as
\begin{align}
    \label{Eq:Tilde-v_r}
    &v_{r}( \nu ) := \int_{ \RR } R_r(x)^2 \nu( \diff x )
    &
    &\forall \nu \in \Pp_2(\RR),
\end{align}
where $R_r$ is the truncating function
\begin{align}
    \label{Eq:R_r}
    R_r(x)
    :&=
    \begin{cases}
        x,	& \Abs{x} \leq r,
        \\
        \frac{r x}{\Abs{x}},	& \Abs{x} > r,
    \end{cases},
\end{align}
for some given parameter $r>0$. Likewise, let $m : \Pp_2(\RR) \to \RR$ be defined as
\begin{align}
    \label{Eq:Tilde-m}
    m(\nu) :&= \int_{\RR} x \nu(\diff x).
\end{align}
It is known from \cite{hernandez-hernandez_coupled_2023} that $m$ and $v_r$ are Lipzchitz functions with respect to the $2$-Wasserstein distance. Consequently, the coefficients
\begin{align}
    \label{Eq:(LQMFG)b,sigma,gamma}
    b_t^\mu( x, a ) :&= b\Parentheses{ x - m(\mu_t) } + a,
    &
    \sigma_t^\mu( x, a ) :&= \sigma,
    &
    \gamma_t^\mu( x, a ) :&= \gamma,
\end{align}
and
\begin{align}
    \label{Eq:(LQMFG)f}
    f_t^\mu( x, a ) 
    :&= 
    \frac{f_1}{2} \Abs{ x - m(\mu_{t}) }^2 
    +
    f \Parentheses{ x - m(\mu_{t}) } a
    +
    \frac{f_2}{2} \Abs{a}^2,
    \\ \label{Eq:(LQMFG)g}
    g^\mu_T( x )
    :&=
    \frac{g}{2} \Abs{ x - m(\mu_{T}) }^2,
\end{align}
are twice differentiable and locally bounded. Note that the $\lambda$-admissible Hamiltionian for the problem,
\begin{align}
    \label{Eq:(LQMFG)H}
    H_t^\nu(x,y,z,u,a)
    :=&
    \frac{f_1}{2} \Abs{ x - m(\nu) }^2 
    +
    f \Parentheses{ x - m(\nu) } a
    +
    \frac{f_2}{2} \Abs{a}^2
    \\ \nonumber
    &+
    b \Parentheses{ x - m(\nu) } y + a y + \sigma z + \gamma u \lambda_t,
\end{align}
is convex in $x$ owing to the fact that $\Diff ^2 H_t^\nu = f_1 > 0$, and it reaches its minimum with respect to $a$ at the point
\begin{align*}
    \widehat{a} = - \frac{f(x - m(\nu)) + y}{f_2}.
\end{align*}
Then, from Remark \ref{Cor:Uncontrolled-jump-diffusion}, there exists a Mean-Field equilibrium $(\widehat{X}, \widehat{\mu})$ to the Linear-Quadratic MFG given by \eqref{Eq:(LQMFG)-J}-\eqref{Eq:LQ-lambda} if, and only if, there exists a solution to the FBSDE 
\begin{align}
    \label{Eq:(LQMFG)-Hat-X}
    \widehat{X}_t 
    =&
    X_0
    +
    \int_0^t \Braces{ 
        b \Parentheses{ \widehat{X}_{s} - \Esp[^1]{\widehat{X}_s} } 
        + 
        \widehat{\alpha}_s 
    } \diff s
    +
    \sigma W_t
    +
    \gamma \tilde{N}^\lambda_t,
    \\
    \label{Eq:(LQMFG)-Hat-Y}
    \widehat{Y}_t
    =&
    g \Parentheses{ \widehat{X}_{T} - \Esp[^1]{\widehat{X}_T} } 
    +
    \int_t^T \Braces{
        f_1 \Parentheses{ \widehat{X}_{s} - \Esp[^1]{\widehat{X}_s} } 
        +
        f \widehat{\alpha}_s 
        +
        b \widehat{Y}_t
    } \diff s
    \\ \nonumber
    &-
    \int_t^T \widehat{Z}_s \diff W_s
    -
    \int_t^T \widehat{U}_{s} \diff \tilde{N}^\lambda_s
\end{align}
where $\widehat{\alpha}$ is the admissible control given by the Pontryagin condition, i.e.
\begin{align}
    \label{Eq:(LQMFG)-Hat-alpha}
    \widehat{\alpha}_t
    \equiv
    \widehat{\alpha}( \Ll^1(\widehat{X})_t, \widehat{X}_t, \widehat{Y}_t )
    :=
    - \frac{ f }{f_2} \Parentheses{ \widehat{X}_{t} - \Esp[^1]{\widehat{X}_t} }
    -
    \frac{ 1 }{f_2} \widehat{Y}_{t}.
\end{align}
However, observe that \eqref{Eq:(LQMFG)-Hat-X}-\eqref{Eq:(LQMFG)-Hat-alpha} form a linear system of FBSDEs. Thus, we need only to verify the $G$-monotonicity condition from \cite[Thm. 7]{hernandez-hernandez_coupled_2023}. That is, we need to verify that there exist constants $\beta_1$, $\beta_2$, $\beta_3 \geq 0$ with $\beta_1 + \beta_2 > 0$ and $\beta_2 + \beta_3 > 0$ such that the inequalities 
\begin{align}
    \label{Eq:(LQMFG)G-monotonicity-I}
    \Parentheses{\delta {A}_t^\mu(x,y,z,u) \cdot \delta (x, y, z, u)}
    &\leq
    - 
    \beta_1 \Abs{ \delta x }^2
    -
    \beta_2 \Braces{
        \Abs{ \delta y }^2
        +
        \Abs{ \delta z }^2
        +
        \Norm{ \delta u }^2_{\lambda_t}
    },
    \\ \label{Eq:(LQMFG)G-Monotonicity-II}
    \Parentheses{\delta g_T^\mu (x) \cdot \delta x}
    &\geq
    \beta_3 \Abs{ \delta x }^2,
\end{align}
hold $\diff \PP \otimes \diff t$-a.e. for all $
    (x,y,z,u), (x',y',z',u')
    \in
    \RR^d \times \RR^n \times \RR^{ d \times k } \times \mathfrak{H}^\lambda
$, where $A^\mu$ is given by
\begin{align*}
    A_t^\mu(x,y,z,u)
    :&=
    \begin{bmatrix}
        -f_1 (x - m(\mu_t)) - f\widehat{\alpha}( \mu_t, x, y ) - by
		\\
		b(x - m(\mu_t)) + \widehat{\alpha}( \mu_t, x, y )
		\\
		\sigma
        \\
        \gamma
    \end{bmatrix}.
\end{align*}
In this case, $\beta_3 = g >0$ due to the convexity of $g_T^\mu$. For $\beta_1$ and $\beta_2$, direct calculations yield 
\begin{align*}
    \Parentheses{\delta {A}_t^\mu(x,y,z,u) \cdot \delta (x, y, z, u)}
    &= 
    - \Parentheses{f_1 - \frac{f^2}{f_2}} \Abs{ \delta x }^2
    -
    \frac{1}{f_2} \Abs{ \delta y }^2
    +
    b( \delta x \cdot \delta y )
    \\
    &\leq
    - \underbrace{\Parentheses{f_1 - \frac{f^2}{f_2} - \frac{\Abs{b}}{2}} }_{=: \beta_1} \Abs{ \delta x }^2
    -
    \underbrace{\Parentheses{\frac{1}{f_2} - \frac{\Abs{b}}{2} } }_{=:\beta_2} \Abs{ \delta y }^2.
\end{align*}
From the hypotheses on the coefficients at \eqref{Eq:(LQFMG)Coefficient-conditions}, inequalities \eqref{Eq:(LQMFG)G-monotonicity-I}-\eqref{Eq:(LQMFG)G-Monotonicity-II} hold; as a result, there exists a strong MFE $(\widehat{\mu},\widehat{X})$,  with $\widehat{\alpha}$ defined in \eqref{Eq:(LQMFG)-Hat-alpha} being the corresponding optimal control. 

\subsection{Portfolio selection with regime switching}
\label{Sec:RS-Porfolio-Selection}

Let $R$ and $S$ be the price of a risk-free bond and a risky asset, respectively. Assume the interest rate $b_0$ and excess return rate $b_1$ offered depend on a grading system $\xi$ that evaluates the current state of the assets on a scale from $1$ to $n$ depending on the current variability of the market,  represented in this case by a given exogenous environment $\mu$.  Mathematically, $R$ and $S$ are defined on an admissible set-up as the solution processes to the (stochastic) differential equations
\begin{align*}
    &\diff R_t = b_0(\xi^\mu_{t-}) R_t \diff t,
    &
    &\diff S_t 
        = \Brackets{ b_0(\xi^\mu_{t-}) + b_1(\xi^\mu_{t-}) } S_t \diff t + \sigma S_t \diff W_t,
\end{align*}
for all $t \in [0,T]$, where 
\begin{itemize}
    \item 
        $R_0 = 1$ and $S_0 \equiv X_0$ such that $\Ll(X_0) = \mu_0$ and $\sigma \neq 0$,

    \item
        $b_0$ and $b_1$ are  maps from $\Ee:= \Braces{ 1, \ldots, n }$ to $(0,\infty)$,

    \item
        $W$ is a Brownian motion independent of $\xi^\mu$,

    \item 
        The regime $\xi^\mu$ depends to the random environment $\mu$ with a finite state space $\Ee$ and transition probabilities
        \begin{align*}
            \Prob{\xi_t^\mu = j \middle\vert \xi^\mu_{t - \delta} = i, \mu_{t-\delta} = \nu}
            =
            \begin{cases}
                \frac{c_r(\nu)}{n} \delta + o(\delta),              & i \neq j,
                \\
                1 - (n-1) \frac{c_r(\nu)}{n} \delta + o(\delta),    & i = j,
            \end{cases}
        \end{align*}
        where $c_r(\nu) = \frac{v_r(\nu) + 1}{r^2 + 1}$.
\end{itemize}
In other words, $\xi^\mu$ is a conditional Markov chain In continuous time and transition semigroup
\begin{align*}
    Q^{(i,j)}(\nu)
    &=
    \frac{c_r(\nu)}{n} \Parentheses{ 1 - n \Ind{ \Braces{i = j} } },
    &
    &\forall \nu \in \Pp_2(\RR^d).
\end{align*}

From our previous discussions on regime-switching diffusions, the chain $\xi^\mu$ can be rewritten as 
\begin{align*}
    \xi_t
    =
    \int_{ \{t\} \times \Ee } e \ \eta^\lambda(\diff s, \diff e),
\end{align*}
where $\eta^\lambda$ is the lifting the admissible common noise process $(N^\lambda, \xi^\lambda)$ with intensity kernel
\begin{align*}
    K^\mu(s,\diff e)
    =
    \sum_{j \in \Ee} Q^{(i,j)}(\mu_{t-})\ \delta_{\{j\}}(\diff e)
    =
    Q^{( i, \cdot )}(\mu_{t-}) \delta_{\Ee}(\diff e).
\end{align*}
Consequently, $\Parentheses{ X_0, \mu, K^\lambda,(N^\lambda,\xi),W }$ forms an admissible set-up.

\begin{remark}
Observe that in this case $\xi^\mu$ can be regarded as a Markov chain that is ``almost'' uniform in their states: the transition probabilities $c_r(\mu_t)/n$ converge to $1/n$ as the variance of the environment increases. Moreover, it is possible to write the intensity candidate $\lambda$  explicitly: 
\begin{align*}
    \lambda_t
    =
    \int_{[0,t) \times \RR} \psi( r, \mu_{s-}, \xi^\mu_{s-} ) N(\diff s, \diff r),
\end{align*}
where $N$ is an independent Poisson random measure on $\RR_+ \times \RR_+$ and $\psi$ is the integrable, $\FF^{\mu,N}$-predictable function
\begin{align*}
    \psi(r,\nu,i)
    :&=
    \sum_{j^+ = 1 }^{n-i} j^+ \Ind{ \Gamma^{+}(\nu,i,j^+) }(r)
    -
    \sum_{j^- = 1 }^{i-1} j^- \Ind{ \Gamma^{-}(\nu,i,j^-) }(r),
\end{align*}
for all $(r,\nu,i) \in \RR_+ \times \Pp_2(\RR^d) \times \Ee$, and $\Gamma^\pm$ are defined as
\begin{align*}
    \Gamma^+(\nu,i,j)
    :=&
    c_r(\nu) \times \left[ (i - 1) + \frac{j-1}{n}, (i - 1) + \frac{j}{n} \right),
    \\
    \Gamma^-(\nu,i,j)
    :=&
    c_r(\nu) \times \left[ (i - 1) - \frac{j}{n}, (i - 1) - \frac{j-1}{n} \right).
\end{align*}
\end{remark}

Let $A = [a_*,a^*] \subset \RR$, $a_* < a^*$, be the amount of money available for an investor to put into the risky asset along the time interval $[0,T]$, and let $\Aa$ be the corresponding set of admissible strategies. Then, for each admissible  strategy $\alpha \in \Aa$, we define the  investment cost (i.e. the 
payoff functional) as the functional
\begin{align}
    \label{Eq:(FBSDE-Toy-regime)J}
    J^\mu(X,\alpha)
    :=&
    \Esp{ 
        \int_0^T \Parentheses{ X_{t-}^2 - m(\mu_{t-})^2 } \frac{f(\xi^\mu_{t-})}{T} \diff t
        +
        \Parentheses{ X_{T-}^2 - m(\mu_{T-})^2 } g(\xi^\mu_{T-})
    },
\end{align}
with $f,g:\Ee \to (0,\infty)$, and  the wealth processes $X \equiv X^\alpha$ follows the controlled SDE
\begin{align}
    \label{Eq:(FBSDE-Toy-regime)X}
    X_t 
    &= 
    X_0
    + 
    \int_0^t 
        \Braces{ b_0\Parentheses{ \xi^\mu_{s-} } X_s + b_1\Parentheses{ \xi^\mu_{s-} } \alpha_s }
    \diff s
    +
    \sigma \int_0^t \alpha_s \diff W_s,
\end{align}
which in turn is well defined due to the compactness of $A$.


Since $b_0,b_1$ and $g$ depend on the environment only through the regime $\xi$, measurability 
and Lipschitz conditions hold for the coefficients
\begin{align}
    \label{Eq:(FBSDE-Toy-regime)b,sigma,gamma}
    b_t^\mu( x, a) :&= b_0\Parentheses{ \xi_{t-}^\mu }x + b_1\Parentheses{ \xi_{t-}^\mu }a,
    &
    \sigma_t^\mu( x, a ) :&= \sigma a,
    &
    \gamma_t^\mu( x, a ) :&= 0,
\end{align}
and
\begin{align}
    \label{Eq:(FBSDE-Toy-regime)f,g}
    f_t^\mu( x, a ) :&= ( x^2 - m(\mu_{t-})^2 ) \frac{ f( \xi_{t-}^\mu ) }{T} ,
    &
    g^\mu_T( x ) :&= ( x^2 - m(\mu_{T-})^2 ) g( \xi_T^\mu ).
\end{align}
Consequently, the $\lambda$-admissible Hamiltonian for the problem is given by 
\begin{align}
    \label{Eq:(FBSDE-Toy-regime)H}
    H_t^\mu(x,y,z,a)
    =
    ( x^2 - m(\mu_{t})^2 )\frac{ f( \xi_{t-}^\mu ) }{T} 
    +
    b_0(\xi_{t-}^\mu) xy
    +
    \Parentheses{ b_1(\xi_{t-}^\mu) y + \sigma z }a,
\end{align}
with the corresponding adjoint equations
\begin{align}
    \label{Eq:(FBSDE-Toy-regime)Y}
    Y_t
    =&
    2 g(\xi_T^\mu) X_{T-}
    +
    \int_t^T \Braces{ \frac{2 f(\xi_{s-}^\mu)}{T} X_{s-} + b_0(\xi_{s-}^\mu) Y_{s-} } \diff s
    -
    \int_t^T Z_s \diff W_s
    \\ \nonumber
    &-
    \int_{(t,T] \times \RR_+} U_s(r) \widetilde{\eta}^\lambda( \diff s, \diff r )
    -
    (M_T - M_t),
    \\ \label{Eq:(FBSDE-Toy-regime)Tilde-Y}
    \widetilde{Y}_t
    =&
    2 g(\xi_T^\mu) 
    +
    2 \int_t^T 
        \Braces{ \frac{f(\xi_{s-}^\mu)}{T} + b_0( \xi_{s-}^\mu ) \widetilde{Y}_{s-} } \diff s
    -
    \int_t^T \widetilde{Z}_s \diff W_s
    \\ \nonumber
    &-
    \int_{(t,T] \times \RR_+} \widetilde{U}_s(r) \widetilde{\eta}^\lambda( \diff s, \diff r )
    -
    (\widetilde{M}_T - \widetilde{M}_t).
\end{align}

Since the Hamiltonian $H$ is linear on $a$, the minimization of \eqref{Eq:(FBSDE-Toy-regime)H} occurs at the extremal points of $A$. Then, the candidate for being an optimal control has the form 
\begin{align}
    \label{Eq:(FBSDE-Toy-regime)hat-alpha}
    \widehat{\alpha}_t
    :=
    \widehat{\alpha}( Y_{t-}, Z_{t-} )
    =
    \begin{cases}
        a^*,    &\mbox{if } b_1(\xi_{t-}^\mu) Y_{t-} + \sigma Z_{t-} \leq 0,
        \\
        a_*,    &\mbox{if } b_1(\xi_{t-}^\mu) Y_{t-} + \sigma Z_{t-} > 0,
    \end{cases}
\end{align}
and the corresponding $\widehat{\alpha}$-controlled Hamiltonian system \eqref{Eq:(FBSDE-Toy-regime)X}-\eqref{Eq:(FBSDE-Toy-regime)Y}-\eqref{Eq:(FBSDE-Toy-regime)Tilde-Y} becomes a Linear FBSDE. 

Lastly, observe that both $g_T^\mu$ in \eqref{Eq:(FBSDE-Toy-regime)f,g} and the mapping 
\begin{align*}
    x 
    \longmapsto 
    b_0(\xi_{t-}^\mu) Y_{t-} x
    +
    \Parentheses{ b_1(\xi_{t-}^\mu) Y_{t-} + \sigma Z_{t-} }\widehat{\alpha}_t
\end{align*}
are convex $\diff t \otimes \diff \PP$-a.a.; thus, in order to 
apply Theorem \ref{Thm:Characterisation-of-Mean-Field-equilibrium}, we only need to verify the $G$-monotonicity conditions \eqref{Eq:(LQMFG)G-monotonicity-I}-\eqref{Eq:(LQMFG)G-Monotonicity-II}. In this case, we have that
\begin{align*}
    A_t^\mu(x,y,z,u)
    :&=
    \begin{bmatrix}
        -\frac{2f(\xi_t^\mu)}{T}x-b_0 y
		\\
		b_0\Parentheses{ \xi_t^\mu }x + b_1\Parentheses{ \xi_t^\mu } \widehat{\alpha}(y,z)
		\\
		\sigma \widehat{\alpha}(y,z)
        \\
        0
    \end{bmatrix},
\end{align*}
and therefore
\begin{align*}
    \Big(
        \delta A_t^\mu( x, y, z, u )
		\cdot 
		\delta ( x, y, z, u )
    \Big)
    &= 
    - \frac{2 f(\xi_{t-}^\mu)}{T} \Abs{ \delta x }^2
    +
    \Parentheses{ b_1(\xi_{t-}^\mu) \delta y + \sigma \delta z } \delta \widehat{\alpha}(y,z).
\end{align*}
However, due to equation \eqref{Eq:(FBSDE-Toy-regime)hat-alpha},
\begin{align*}
    \Parentheses{ b_1(i) \delta y + \sigma \delta z } \delta \widehat{\alpha}(y,z)
    \leq 
    0
\end{align*}
for all $y, z \in \RR$ on the event $\Braces{ \xi_t^\mu = i }$ for all $i \in \Ee$. Therefore, by taking $\beta_1 = T^{-1} \min_{i \in \Ee} 2 f(i)  > 0$, $\beta_2 = 0$ and $\beta_3 = \min_{i \in \Ee} \Braces{2 g(i)} >0$, $G$-monotinicity holds.

\end{document}